\documentclass[draft]{amsart}
\usepackage{amssymb}
\def\C{\mathbb{C}}

\def\R{\mathbb{R}}

\def\K{\mathcal K}

\def\B{\mathcal B}

\def\a{\alpha}

\def\beqa{\begin{eqnarray*}}
\def\eeqa{\end{eqnarray*}}

 \newtheorem{thm}{Theorem}
 \newtheorem{cor}[thm]{Corollary}
 \newtheorem{lem}[thm]{Lemma}
 \newtheorem{prop}[thm]{Proposition}
 \newtheorem{defn}[thm]{Definition}
 \newtheorem{rem}[thm]{Remark}
 \newtheorem{ex}[thm]{Example}
\newcommand{\be}{\begin{equation}}
\newcommand{\ee}{\end{equation}}
\newcommand{\bea}{\begin{eqnarray}}

\newcommand{\eea}{\end{eqnarray}}
\newcommand{\Bea}{\begin{eqnarray*}}
\newcommand{\Eea}{\end{eqnarray*}}

\newcounter{cnt1}
\newcounter{cnt2}
\newcounter{cnt3}
\newcommand{\blr}{\begin{list}{$($\roman{cnt1}$)$}
 {\usecounter{cnt1} \setlength{\topsep}{0pt}
 \setlength{\itemsep}{0pt}}}
\newcommand{\bla}{\begin{list}{$($\alph{cnt2}$)$}
 {\usecounter{cnt2} \setlength{\topsep}{0pt}
 \setlength{\itemsep}{0pt}}}
\newcommand{\bln}{\begin{list}{$($\arabic{cnt3}$)$}
 {\usecounter{cnt3} \setlength{\topsep}{0pt}
 \setlength{\itemsep}{0pt}}}
\newcommand{\el}{\end{list}}

\date{}

\begin{document}

\author[Bhat, B. V. R.]{B. V. Rajarama Bhat}
\author[Lindsay, J. M.]{Martin Lindsay}
\author[Mukherjee, M.]{Mithun  Mukherjee}
\title{ Additive units of product system}
\address[Bhat, B. V. R.]{Indian Statistical Institute, R. V. College Post, Bangalore-560059, India}
\email{bhat@isibang.ac.in}

\address[Lindsay, J. M.]{Lancaster University, Lancaster- LA1 4YW, UK}
\email{jmartinlindsay@googlemail.com}

\address[Mukherjee, M.]{IISER Kolkata, Mohanpur-741 252, India}
\email{mithun.mukherjee@iiserkol.ac.in}

\maketitle
\vskip 10pt

\begin{center}
{\bf Abstract}

\end{center}

\vskip 4pt

\begin{abstract}

We introduce the notion of additive units and roots of a unit in a spatial product system. The set of all roots of any unit forms a Hilbert space and its dimension is the same as the index of the product system. We show that a unit and all of its roots generate the type I part of the product system. Using properties of roots, we also provide an alternative proof of the Powers' problem that the cocycle conjugacy class of Powers sum is independent of the choice of intertwining isometries. In the last section, we introduce the notion of cluster of a product subsystem and establish its connection with random sets in the sense of Tsirelson (\cite{Tsi-three-questions}) and Liebscher (\cite{Lie-random}).

\end{abstract}

\section{Introduction}

A fundamental goal of quantum dynamics is the classification of semigroups of unital $*$-endomorphisms of the algebra of all bounded operators on a separable Hilbert space up to cocycle conjugacy. Associated with every such `$E_0$-semigroup', is a (tensor) product system of Hilbert spaces (\cite{Arv-continuous}). This translates the problem of classification of $E_0$-semigroups up to cocycle conjugacy  into the problem of classification of the product systems up to isomorphism. A product system is a measurable family of separable Hilbert spaces $(\mathcal E_s)_{s>0}$ with associative identification $\mathcal E_{s+t} \simeq \mathcal E_s\otimes \mathcal E_t$ through unitaries. A unit is a measurable section of non-zero vectors $(u_s)_{s>0},$ $u_s \in \mathcal E_s$ which factorises: $u_{s+t}=u_s\otimes u_t,$ $s,t>0.$ Depending on the existence of units, product systems are classified into three categories. A product system is said to be of type I if units exist and they `generate' the product system. A product system is said to
be of type II  if it has a unit but they fail to `generate' the product system. Product systems having units are also known as spatial product systems. A product system is said to be of type III or non-spatial if it does not have any unit. Spatial product systems have an index. The index is a complete invariant for type I product systems and each is cocycle conjugate to a CCR flow (\cite{Arv-continuous-III}). There is an operation of tensoring on the category of product systems. The index is additive under the tensor product of spatial product systems. Product systems of type II and type III exist in abundance but their classification theory is far from complete. It was shown that there are uncountably many cocycle conjugacy classes of type II and type III product systems (\cite{Pow-nonspatial},\cite{Pow-new-example},\cite{Tsi-non-isomorphic},\cite{Tsi-coloured}) but we still lack good invariants to distinguish them.

Tsirelson (\cite{Tsi-three-questions},\cite{Tsi-Warren's-noise}) established interesting new
examples of type II product systems coming from measure types of random sets
or generalized random (Gaussian) processes. Liebscher, (\cite{Lie-random}) then
made a systematic study of measure types of random sets. Given a pair of product
systems, one contained in the other, one associates a measure type of
random (closed) sets of the interval $[0,1].$ These measure types are stationary
and factorizing over disjoint intervals. The corresponding measure type is an invariant of the product system. See \cite{Lie-random} for more details.

Contractive semigroups of completely positive maps are known as quantum dynamical semigroups. The dilation theory of quantum dynamical semigroups (\cite{Bha-minimal}) reveals a new approach to understand $E_0$-semigroups. Every unital quantum dynamical semigroup dilates to an $E_0$-semigroup and the minimal dilation is unique up to conjugacy.

Similarly, $E_0$ semigroups on general $C^*$ algebras or von Neumann algebras correspond to product systems of Hilbert modules, (\cite{MuS-Markov},\cite{ShS-subproduct},\cite{Ske-classification}). Much of the theory of product system of Hilbert spaces and the theory of $E_0$-semigroups acting on $\B(H)$ can be carried through also for the product systems of Hilbert modules and $E_0$ semigroups acting on $\mathcal B^a(E),$ the algebra of all adjointable operators on a Hilbert module. However there is no natural tensor product operation on the category of product systems of Hilbert modules. Skeide (\cite{Ske-index}) overcame this by introducing the spatial product of spatial product systems of Hilbert modules in which the reference units (normalized) are identified and under which the index of the spatial product system of Hilbert module is additive. Restricting to the case of spatial product systems of Hilbert spaces, we have another operations on the category of spatial product systems. Suppose $\mathcal E$ and $\mathcal F$ are two spatial product systems with normalized units $u$ and $v$ respectively. The spatial product can be identified with the product
subsystem of the tensor product, generated by the two subsystems $\mathcal E\otimes v$ and $u\otimes \mathcal F.$ This raises the question whether the spatial product is the tensor product or not. Powers (\cite{Pow-addition}) answered this in the negative sense by solving the seemingly different but equivalent following problem:

Suppose $\phi =\{ \phi _t: t\geq 0\} $
and $\psi = \{ \psi _t: t\geq 0\} $ are two $E_0$ semigroups on
$\mathcal {B} ( {H})$ and $\mathcal {B}( {K})$ respectively and
$U=\{ U_t: t\geq 0\}$ and $V=\{ V_t: t\geq 0\}$ are two strongly
continuous semigroups of isometries which intertwine $\phi _t$ ($\phi_t(A)U_t=U_tA,$ $\forall ~A \in \B(H),t\geq 0$) and
$\psi _t$ respectively. Note that the intertwining isometries of $E_0$-semigroups correspond bijectively to the normalized units of the associated product systems. Consider the CP semigroup (Powers sum) $\tau _t$ on
$\mathcal {B}( {H}\oplus {K})$ defined by
$$\tau _t\left(
 \begin{array}{cc}
 X & Y \\
 Z & W \\
 \end{array}
 \right)=\left(
 \begin{array}{cc}
 \phi _t(X) & U_tYV^*_t \\
 V_tZU^*_t & \psi _t(W) \\
 \end{array}
 \right).$$
How is the product system of the minimal dilation (in the sense of
\cite{BhP-Markov},\cite{Bha-minimal}) of $\tau  $ related to the product systems of $\phi $ and $\psi $? Skeide (\cite{Ske-commutant}) identified the product system as a spatial product
through normalized units. The definition of Powers' sum easily extends to CP semigroups and the
product system of Powers' sum in that case also is the spatial
product of the product systems of its summands (\cite{BLS-subsystems},\cite{Ske-CPD}). Motivated by this problem and its straightforward generalization to more general `corner', amalgamated product (see Section 2)  through general contractive morphism of two product systems (not necessarily spatial) was introduced in \cite{BhM-inclusion} which
generalizes the spatial product. The spatial product may
be viewed as an amalgamated product through the contractive morphism
defined through normalized units. This answers Powers' problem for
the Powers' sum obtained from not necessarily isometric intertwining
semigroups.

The structure of the spatial product, a priori depends on the choice of the reference units in their respective factors. In fact, Tsirelson (\cite{Tsi-automorphism}) showed that the group of all automorphisms of a product system may not act transitively on the set of all units. It raises another question whether the isomorphism class of the spatial product depends on the choice of the reference units. Equivalently, whether the cocycle conjugacy class of the minimal dilation of Powers sum depends on the choice of the intertwining isometries. This was answered in the negative sense in \cite{BLMS-intrinsic}. See also \cite{BLS-spatial}.

In this paper, we start with a brief overview of the theory of inclusion systems and amalgamated products to make the readers familiar with these notions which we use repeatedly. Readers are referred to \cite{BhM-inclusion}, \cite{Muk-index} for more details. In Section 3, we introduce the notion of additive units and roots of a unit in a spatial product system. Additive units are measurable sections of product system which are `additive with respect to a given unit'. Roots are the special additive units such that for each $t>0,$ the sections are orthogonal to the unit. The set of all additive units forms a Hilbert space and the set of all roots is a subspace of co-dimension one. We compute all the roots of the vacuum
unit in CCR flows ($\Gamma_{sym}(L^2[0,t],K) $). They are given by the set of all $c\chi|_{t]},$ $c\in K$ almost surely. From this, we establish that a unit and all of its roots `generate' the type I part of the product system and the dimension of the Hilbert space of the set of all roots of a unit is the same for every unit and coincides with the index of the product system. We also generalize the notion of additive units and roots of a unit on the level of inclusion systems (see Section 2). We show that the set of all additive units of a unit in an inclusion system are
in a bijective correspondence with the set of all additive units of the `lifted' unit in the generated algebraic product system. The behaviour of the roots under amalgamated product is also studied. Using the properties of roots, we have an alternating proof of the fact that the Powers sum is independent of the choice of the intertwining isometries or equivalently that the isomorphism class of the amalgamated product through normalized units is independent of the choice of the units (see Section 4). In fact, we have an improvement of this result which says that the isomorphism class of the amalgamated product through strictly contractive units is also independent of the choice of the units. This fact will be explained elsewhere (\cite{Muk-contractive}).

In Section 5, given any product subsystem $\mathcal F$ of a product system $\mathcal E,$ we construct an intermediate subsystem called the cluster subsystem of $\mathcal F.$ A product subsystem corresponds to an `adapted' family of commutative projections satisfying some relation. The commutative von Neumann algebra generated by them is uniquely determined by a measure type of random closed sets of the interval $[0,1].$  The distribution of the random mapping which sends a closed set to its limit points is the measure type of the cluster system of the original product subsystem. In a special case, the measure type corresponding to a single unit and the measure type corresponding to the type I part, both share the same relation. See Proposition 3.33, Chapter 3, \cite{Lie-random}. Liebscher's proofs of those facts use heavy machinery from measure theory of random sets and the direct integral construction.  Here we explicitly construct the
cluster subsystem without involving any heavy machinery. We show that the measure type corresponding to the subsystem and the measure type of its cluster are related by the above random mapping. Without using any random sets theory, we also compute that the cluster of the subsystem generated by a single unit in a spatial product system is the type I part of the product system.

\section{Inclusion system and amalgamation}
An inclusion system is a parametrized family of Hilbert spaces exactly like product system but the connecting maps are now only isometries. These objects seem to be ubiquitous in the field of product system. They are the recurrent theme of studying quantum dynamics, in particular CP semigroups. (See \cite{BhS-tensor},\cite{MuS-Markov},\cite{Mar-CP},\cite{ShS-subproduct},\cite{BhM-inclusion}). Even while associating product systems to CP semigroups what one gets first are inclusion systems, and then an inductive limit procedure gives
product systems (\cite{BhS-tensor},\cite{BhM-inclusion}). The notion of inclusion systems is introduced in \cite{BhM-inclusion}. It was also introduced by Shalit and Sholel (\cite{ShS-subproduct}) under the name subproduct system. The following definition is taken from \cite{BhM-inclusion}.
  \begin{defn}
 An inclusion System $(E, \beta )$ is  a family of Hilbert spaces $E=
 \{ E_t,t\in (0, \infty )\}$ together with isometries
 $\beta_{s,t}$:$E_{s+t}\rightarrow E_s\otimes E_t$, for   $s,t \in
 (0, \infty )$, such that $\forall$ $r,s,t \in (0, \infty),$ ~
 $(\beta_{r,s}\otimes 1_{E_t})\beta_{r+s,t}=(1_{E_r}\otimes
 \beta_{s,t})\beta_{r,s+t}.$ It is said to be an algebraic product system if
 further every $\beta _{s, t}$ is a unitary.
 \end{defn}

 \begin{defn}
 Suppose $(E,\beta)$ is an inclusion system. Then a family  $F=(F_t)_{t>0}$ of closed subspaces, $F_t \subset E_t$ is said to be an inclusion subsystem of
 $(E,\beta)$ if
 $\beta_{s,t}|_{F_{s+t}}(F_{s+t})\subset F_s\otimes F_t$ for every
 $s,t>0.$
 \end{defn}

  For each $t \in \R_+,$
 we set $$J_t=\{(t_1,t_2,\ldots , t_n): t_i>0, \sum _{i=1}^n t_i=t, n\geq
 1\}.$$ For $\textbf{s}=(s_1,s_2,\ldots , s_m) \in J_s$, and
 $\textbf{t}=(t_1,t_2, \ldots , t_n) \in J_t$ we define
 $\textbf{s} \smile \textbf{t}:=(s_1,s_2,\ldots ,s_m,t_1,t_2,
 \ldots , t_n) \in J_{s+t}$. Now fix $t \in \R_+$. On $J_t,$ define a
 partial order $\textbf{t}\geq\textbf{s}=(s_1,s_2,\ldots, s_m)$
 if for each $i,$ $(1\leq i\leq m)$ there exists (unique)
 $\textbf{s}_i \in J_{s_i}$ such that
 $\textbf{t}=\textbf{s}_1\smile\textbf{s}_2\smile \cdots
 \smile\textbf{s}_m$. The order relation $\geq$ makes $J_t$ a directed set.

  Suppose $(E,\beta)$ is an
 inclusion system. For $\textbf{s}=(s_1,\cdots,s_n)\in J_t,$ we set $E_\textbf{s}=E_{s_1}\otimes\cdots\otimes E_{s_n}.$ For $\textbf{s}=(s_1,\cdots,s_n)\leq \textbf{t}=\textbf{s}_1\smile\cdots\smile\textbf{s}_n\in J_t,$ define $\beta_{\textbf{s},\textbf{t}}:E_\textbf{s}\rightarrow E_\textbf{t}$ by $\beta_{\textbf{s},\textbf{t}}=\beta_{s_1,\textbf{s}_1}\otimes\cdots\otimes\beta_{s_n,\textbf{s}_n},$ where $\beta_{s,s}=I_{E_s}$ and for $\textbf{s}=(s_1,\cdots,s_n)\in J_s,$ inductively define $$\beta_{s,\textbf{s}}=(I\otimes \beta_{s_{n-1},s_n})\cdots(I\otimes \beta_{s_2,s_3+\cdots+s_n})\beta_{s_1,s_2+s_3+\cdots+s_n}.$$  Proof of the following theorem can be found in Theorem 5, \cite{BhM-inclusion}. \begin{thm}\label{inductive }
 Suppose $(E, \beta )$ is an inclusion system. Let $\mathcal
 E_t=\mbox{indlim}_{J_t} E_\textbf{s}$ be the inductive limit of
 $E_\textbf{s}$ over $J_t$ for $t>0.$ Then $\mathcal E= \{\mathcal
 E_t : t>0\}$ has the structure of an algebraic product system.
 \end{thm} Let $(\mathcal E,B)$ be the generated algebraic product
 system of the inclusion system $(E,\beta).$ Note that the unitary map $B_{s,t}$ goes from $\mathcal E_{s+t}$ to $\mathcal E_s \otimes \mathcal E_t$ for every $s,t>0.$ In other words, algebraic product systems are inclusion systems with all the linking maps are unitaries. Observe that any product system is an algebraic product system but the converse may not be true. The multiplication operation of a product system $\mathcal E$ gives rise to the unitary maps which goes from $\mathcal E_s\otimes \mathcal E_t$ to $\mathcal E_{s+t}$ for every $s,t>0.$ Ad-joints of these unitary maps obviously associative and makes it into an algebraic product system. Therefore we can assume that a product system is a special algebraic product system. Though the linking maps implement `co-product' rather than `product' but abusing of terminology, we call it an algebraic product system. Nevertheless, we can talk about an inclusion subsystem of a product system. The following important fact that an inclusion subsystem in a product
system generates a product subsystem is used throughout without reference. For the proof, see Lemma \ref{subinclusion}, Appendix A. The following definition is taken from \cite{BhM-inclusion}.

 \begin{defn}
  Let $(E, \beta )$ be an inclusion system. Let $u=\{ u_t: t>0\}$ be a
  family of vectors such that $(1)$ for all $t>0,$ $u_t\in E_t$ $(2)$
  there is a $k\in \R,$ such that $\|u_t\|\leq \exp(tk),$ for all $t>0.$
  and  $(3)$ $u_t\neq 0$ for some $t>0.$ Then $u$ is said to be a unit if
$$u_{s+t}= \beta _{s, t}^*(u_s\otimes u_t) ~~\forall s, t>0.$$
  \end{defn}

Let $i_t: E_t\rightarrow \mathcal E_t$ be the canonical embedding.

 \begin{thm}\label{unit}
 Let $(E , \beta)$ be an inclusion system and let $(\mathcal E , B)$
 be the algebraic product system generated by it. Then the map $i^*$ provides a bijection between
the set of all units of $(\mathcal E, B)$ and the set of all units of $(E, \beta )$ by letting it acts point-wise on units.
 \end{thm}

 For the proof, readers are referred to Theorem 10, \cite{BhM-inclusion}.

 Fix a unit $u$ of $(E,\beta).$ Then by the above theorem there is a unique unit $\hat{u}$ in $(\mathcal E,B)$ such that for every $t>0,$ $i^*_t(\hat{u}_t)=u_t.$ We say $\hat{u}$ as the `lift' of $u.$ Note that if $u$ is normalized, then $\hat{u}$ is also normalized.

 \textbf{Amalgamation}

The amalgamated product of two product systems over a contractive morphism is introduced in \cite{BhM-inclusion}. The index of the amalgamated product over general contractive morphism is computed in \cite{Muk-index}. The following theorem characterizes the amalgamated product. See Theorem 2.7, \cite{Muk-index}.

 \begin{thm}\label{universal}
 Suppose $(\mathcal E,W^{\mathcal E})$ and $(\mathcal F,W^{\mathcal
 F})$ are two product systems and let $C:(\mathcal F,W^{\mathcal
 F})\rightarrow (\mathcal E,W^{\mathcal E})$ be a contractive
 morphism. Then there exist an algebraic product system $(\mathcal G,W^{\mathcal G})$ and
 isometric product system morphisms $I:\mathcal
 E\rightarrow \mathcal G$ and $J:\mathcal F\rightarrow \mathcal G$
 such that the following holds:

 (i) $\langle I_s(x),J_s(y)\rangle=\langle x,C_sy\rangle$ for all
 $x\in \mathcal E_s$ and $y\in \mathcal F_s.$

 (ii) $\mathcal G=I(\mathcal E)\bigvee J(\mathcal F).$

 \end{thm}
$\mathcal G$ is said to be the amalgamated product of $\mathcal E$ and $\mathcal F$ over the contractive morphism $C$ and denoted by $\mathcal G=\mathcal E\otimes_C\mathcal F.$ For the details of construction, we refer to Section 3, \cite{BhM-inclusion}.

\section{Additive units}


Suppose $\mathcal E$ is a product system. The multiplication operation in $\mathcal E$ is as follows: For $s,t>0,$ $a \in \mathcal E_s,$ $b \in \mathcal E_t,$ we have $a\cdot b \in \mathcal E_{s+t}$ and $\mathcal E_{s+t}=\overline{\mbox{span}}~\mathcal E_s\cdot\mathcal E_t.$ Also for $a,a^\prime \in \mathcal E_s,$ $b,b^\prime \in \mathcal E_t,$ we have $$\langle a\cdot a^\prime,b\cdot b^\prime \rangle_{\mathcal E_{s+t}} = \langle a,b\rangle_{\mathcal E_s} \langle b,b^\prime \rangle_{\mathcal E_t}.$$ In this section, we abbreviate the multiplication $a\cdot b$ as $ab.$

\begin{defn} Let $\mathcal E$ be a spatial product system and let $u=(u_t)_{t>0}$ be a unit of $\mathcal E.$ A
measurable section $(a_t)_{t>0}$ of $\mathcal E$ is said to be an
additive unit of $u$ if for all $s,t>0,$
$$a_{s+t}=a_su_t+u_sa_t.$$
\end{defn}

\begin{defn}An additive unit $a=(a_t)_{t>0}$ of a unit $u=(u_t)_{t>0}$ is said to be
a root if $\langle a_t,u_t\rangle=0$ for all $t>0.$
\end{defn}

\begin{rem}
It is clear that the set of all additive units of a given unit forms a vector space under point wise addition and point wise scalar multiplication. The set of all roots forms a vector subspace of it. Indeed if $a=(a_s)_{s>0}$ and $b=(b_s)_{s>0}$ are two additive units(roots) of a unit $u,$ then clearly $\lambda a:=(\lambda a_s)_{s>0}$ and $(a+b):=(a_s+b_s)_{s>0}$ are additive units(roots) of $u=(u_s)_{s>0}.$ Also note that, if $a$ is an additive unit(roots) of $u,$ then $(a^\prime)_{s>0}$ which is defined by $a^\prime_s=exp(\lambda s)a_s,$ is an additive unit(roots) of $(u^{\prime}_s)_{s>0},$ where $u^\prime_s= exp(\lambda s)u_s.$ In other words, the additive units of a unit are completely determined by the additive units of the normalized unit.
\end{rem}

\begin{ex}
Let $u=(u_s)_{s>0}$ be a unit in a product system $\mathcal E.$ Then the measurable section $b=(b_s)_{s>0}$ given by $b_s=\lambda su_s,$ for some $\lambda \in \mathbb C,$ $s>0,$ is an additive unit of the unit $u.$ We call them as the trivial additive units of the unit $u.$
\end{ex}

Let $a$ be an additive unit of a unit $u.$ For $s>0,$ consider the measurable function $$f:\mathbb R_+\rightarrow \mathbb C$$ given by $$f(s)=\langle u_s,a_s\rangle \|u_s\|^{-2}.$$ Then a simple computation shows that $f(s+t)=f(s)+f(t),$ $s,t>0.$ This implies $f(s)=sf(1).$ Decomposing $a_s=b_s+b^\prime_s,$ where $$b_s=\langle u_s,a_s\rangle \|u_s\|^{-2}u_s$$  and $$b^\prime_s=a_s-(\langle u_s,a_s\rangle\|u_s\|^{-2} u_s),$$ we find that $b_s=(\lambda su_s)_{s>0}$ for some $\lambda \in \mathbb C$ and $b^\prime$ is a root of $u.$ In other words, every additive unit decomposes uniquely as a trivial additive unit and a root. From the remark, we may assume without loss of generality that our unit $u$ is normalized, i.e. $\|u_s\|=1,$ for every $s>0.$ Let $a$ and $b$ be two roots of the normalized unit $u.$  Then a similar computation shows that $$\langle a_s,b_s\rangle=s \langle a_1,b_1\rangle,~ s>0.$$ Now consider $a,b$ two additive units of $u.$ Then we can decompose $$a_s=c_s+c^\prime_s~,~b_s=d_s+d^\prime_s~,~
s>0,$$ where $$c_s= s \langle u_1,a_1\rangle u_s~,~d_s= s \langle u_1,b_1 \rangle u_s~,s>0,$$ and $c^\prime,d^\prime $ are roots of $u$ with $$\langle c^\prime_s,d^\prime_s \rangle= s \langle c^\prime_1,d^\prime_1 \rangle.$$ Now $\langle c^\prime_1,d^\prime_1\rangle=\langle(a_1-\langle u_1,a_1\rangle u_1),(b_1-\langle u_1,b_1\rangle u_1)\rangle=\langle a_1,b_1\rangle-\langle a_1,u_1\rangle \langle u_1,b_1 \rangle.$ From this, a simple computation shows $\langle a_s,b_s \rangle = s^2\langle a_1,u_1 \rangle \langle u_1,b_1\rangle + s \langle a_1,b_1 \rangle-s \langle a_1,u_1 \rangle \langle u_1,b_1\rangle.$ In other words,\bea \label{roots} \langle a_s,b_s\rangle=\langle \theta_s a_1,\theta_sb_1\rangle,\eea where $\theta_s:\mathcal E_1\rightarrow \mathcal E_1$ is given by $\theta_s=[sI+(s^2-s)|u_1><u_1|]^\frac{1}{2}.$
\begin{prop}
Let $u$ be a normalized unit of a product system $\mathcal E$. Then the set of all additive units of $u$ forms a Hilbert space under the inner product $\langle a,b \rangle =:\langle a_1,b_1 \rangle_{\mathcal E_1}$ and the set of all roots of $u$ is a closed subspace of co-dimension one.
\end{prop}

Proof:  Let us denote by $A^\mathcal E_u$ and $R^\mathcal E_u$ be the vector spaces of all additive units and roots of $u$ respectively. For $a,b \in A^\mathcal E_u,$ define an inner product on $A^\mathcal E_u$ by $\langle a,b\rangle=\langle a_1,b_1\rangle.$  Let $\{a^n\}_{n\geq 1}$ be a Cauchy sequence. i.e. $\|a^n-a^m\|\rightarrow 0$ as $m,n\rightarrow \infty.$ Now from Equation \ref{roots}, we get for each $s>0$ $\|a^n_s-a^m_s\|=\|\theta_s(a^n_1-a^m_1)\|\leq\|\theta_s\|\|a^n_1-a^m_1\|=\|\theta_s\|\a^n-a^m\|\rightarrow 0.$ Let for each $s>0,$ $a_s=\mbox{Lim}_{n\rightarrow \infty}a^n_s.$ The section $(a_s)_{s>0}$ is clearly measurable as being point-wise limit of measurable sections. Now we will show that $(a_s)_{s>0}$ is an additive unit of $u.$ Let $\epsilon>0$ be given. For $s,t>0$ choose N such that for $n>N,$ $\|a^n_s-a_s\|\leq \frac{1}{3}\epsilon,$ $\|a^n_t-a_t\|\leq \frac{1}{3}\epsilon $ and $\|a^n_{s+t}-a_{s+t}\|\leq \frac{1}{3}\epsilon.$ Then \allowdisplaybreaks{\begin{align*} & \|a_{s+t}-a_su_t-u_sa_t\| \\ &
\leq  \|a_{s+t}- a^n_{s+t}\|+\|a^n_su_t-a_su_t\|+\|u_sa^n_t-u_sa_t\| \leq \epsilon .\end{align*}}%
So $a\in A^\mathcal E_u$ and $\|a^n-a\|\rightarrow 0.$ This proves that $A^\mathcal E_u$ is complete with respect to the inner product. Other part is trivial. \qed

\begin{prop}\label{vacuum} The set of all roots of the vacuum unit in CCR flow of index $k$ is given by the set $\{c\chi_{t]}:c\in K\}$ almost everywhere.
\end{prop} Proof: It is easy to see that $c\chi_{t]}, c\in K$ are the roots of the vacuum unit. To prove the converse, if $a$ is
a root of the vacuum unit, then in Guichardet's picture described in Appendix B , the following identity is valid
almost everywhere,
$$a_{s+t}(\sigma)= \left\{
                              \begin{array}{ccc}
                                a_s(\sigma \cap [t,s+t]-t) & \mbox{if}~{~} \sigma\cap[0,t]=\emptyset  \\
                                a_t(\sigma \cap [0,t]) & \mbox{if}~{~} \sigma\cap[t,s+t]=\emptyset  & ~,{~}\sigma\in \Delta(s+t)\\
                                0 & \mbox{otherwise}. \\
                              \end{array}
                            \right.$$
Fix $s\in \mathbb R_+.$ Let $k$ be any natural number. Denote by $[a,b]^\prime$ the complement of $[a,b]$ in $[0,s].$ Then we have the identity,
almost everywhere ,$$ a_s(\sigma)=\left\{
                                                                           \begin{array}{ccc}
                                                                             a_\frac{s}{k}(\sigma\cap[\frac{(k-1)s}{k},s]-\frac{(k-1)s}{k}) & \mbox{if}~{~} \sigma\cap[\frac{(k-1)s}{k},s]^\prime=\emptyset \\
                                                                             a_\frac{s}{k}(\sigma\cap[\frac{(k-2)s}{k},\frac{(k-1)s}{k}]-\frac{(k-2)s}{k}) & \mbox{if}~{~} \sigma\cap [\frac{(k-2)s}{k},\frac{(k-1)s}{k}]^\prime=\emptyset \\
                                                                             . & . & ~,{~}~ \sigma\in \Delta(s) \\
                                                                             . & . \\
                                                                             a_\frac{s}{k}(\sigma\cap[0,\frac{s}{k}]) & \mbox{if}~{~} \sigma\cap[0,\frac{s}{k}]^\prime=\emptyset  \\
                                                                             0 &
                                                                             \mbox{otherwise}.
                                                                           \end{array}
                                                                         \right.
 $$

                            Suppose that $\#~\sigma=n,$ then the subset of $\Delta_n(s),$ where $a_s$ is non zero except on a set of measure zero, is contained in
$$\cup_{i=0}^{k-1}[\Delta_n(s/k)+is/k],~{~} \mbox{for}~\mbox{all}~k=1,2,\cdots.$$ The Lebesgue  measure of the set $\cup_{i=0}^{k-1}[\Delta_n(s/k)+is/k]$ is $s^n/n!k^{n-1}.$ So the Lebesgue measure of the set $\cap_{k\geq 1}\cup_{i=0}^{k-1}[\Delta_n(s/k)+is/k]$ is zero for $n\geq 2.$ It follows that $a_s$ vanishes on $\Delta_n(s),$ for $n\geq 2.$
As it is a root, it is orthogonal to the vacuum unit,  we conclude
that, $a_s$ is a measurable function in $L^2([0,s],K)$ with the
property, a.e. \bea \label{shift} a_s=a_r+S_ra_{s-r} ,~~\forall r,~
0<r<s.\eea where $S_t$ on $L^2(\mathbb R_+,K)$ defined by
\bea\label{S_t}(S_tf)(s)&=&\left\{
                                                                                                        \begin{array}{cc}
                                                                                                          f(s-t) & \mbox{if}~s\geq t \\
                                                                                                          0 & \mbox{otherwise}. \\
                                                                                                        \end{array}
                                                                                                      \right.\eea For every $x\in K,$  define the measurable function
$A_x:\mathbb R_+\rightarrow \mathbb C,$ by $A_x(s)=\langle
a_s,x\chi|_{s]}\rangle.$ An easy calculation shows that
$A_x(s+t)=A_x(s)+A_x(t).$ Its measurable solution is given by
$A_x(s)=sA_x(1).$ Let us define the linear functional
$f:\K\rightarrow \mathbb C$ by $f(x)=\langle
a_1,x\chi|_{1]}\rangle,$ for $x\in K.$ It is bounded as $\|f\|\leq
\|a_1\|_{L^2}.$ So by Riesz representation theorem there is a unique
$y\in K$ such that $f(x)=\langle y,x\rangle.$  Now for $r\leq s,$
$z\in K,$ \allowdisplaybreaks{\begin{align*} \langle a_s-y\chi|_{s]}, z\chi|_{r]}\rangle &=
\langle a_r,z\chi|_{r]}\rangle-r\langle y,z\rangle\\
&= rA_z(1)-r\langle y,z\rangle\\ &= 0. \end{align*}} %
As the set
$\{z\chi|_{r]}: z\in K, 0\leq r\leq s\}$ is total in $L^2([0,s],K),$
we have $a_s=y\chi|_{s]}.$ \qed

Let us denote by $R^\mathcal E_u,$  the Hilbert space of roots of the unit $u$ in $\mathcal E.$

\begin{thm}
Suppose $(\mathcal E,W)$ is a product system and $u$ is a normalized unit of $\mathcal E.$ Then $\mbox{dim}~R^\mathcal E_u = \mbox{index}~\mathcal E.$
\end{thm}

Proof: First we claim that roots of $u$ are in $\mathcal E^I,$ the type I part of $\mathcal E.$  Given a root $a$ of $u,$ $\|a\|=1,$ set $E_s= \mbox{span}~\{u_s,a_s\}.$ Then it is easy to see that $(E,W|_E)$ is an inclusion system. Let $\Gamma_{sym}(L^2[0,t])$ be the symmetric Fock product system. Define $\phi_s : E_s \rightarrow \Gamma_{sym}(L^2[0,t]$ by $\phi_s(\alpha u_s + \beta a_s)=\alpha \Omega_s + \beta \chi|_{s]}.$ Then $\phi=(\phi_s)_{s>0}$ is an isometric morphism of inclusion system. So the product system generated by $u$ and $a$ is isomorphic to a type I product system. This proves the claim. Any isomorphism of $\mathcal E^I$ to $\Gamma_{sym}(L^2[0,t],K)$ ($\mbox{dim}~K=\mbox{index}~\mathcal E$) sending $u$ to vacuum unit, sends roots to roots. This implies every root of $u$ under this will be mapped to a $(c\chi|_{s]})_{s>0}$ and vice versa. The result now follows. \qed

\begin{cor}\label{root I}
Let $a$ be a root of a unit $u$ in a spatial product system
$(\mathcal E,B).$ Then $a \in \mathcal E^I.$
\end{cor}

\begin{cor}\label{unit-roots}
Suppose $(\mathcal E,B)$ is a spatial product system and $u$ is a
unit. Then the product system generated by the unit $u$ and all roots of
it, is the type I part of $(\mathcal E,B).$
\end{cor}

We shall now define all these notions on the level of inclusion system. We quote the following definition from \cite{BhM-inclusion}.

\begin{defn}
Let $(E,\beta)$ be an inclusion system and let $u$ be a normalized unit of $(E,\beta).$ A section $(a_t)_{t>0}$ of $(E,\beta)$ is said to be an additive unit of the unit $u$ if $$a_{s+t}=\beta^*_{s,t}(a_s\otimes
u_t+u_s\otimes a_t)~~ \mbox{and}~~ \|a_s\|^2\leq k(s+s^2),~s>0 ,~
\mbox{for}~~ \mbox{some}~~ k\geq 0.$$
\end{defn}

\begin{defn}An additive unit $a=(a_t)_{t>0}$ of a unit $u=(u_t)_{t>0}$ is said to be
a root if $\langle a_t,u_t\rangle=0$ for all $t>0.$
\end{defn}

\begin{prop}\label{root}
Let $(E , \beta)$ be an inclusion system and let $(\mathcal E , B)$
be the algebraic product system generated by it. Then  $i^*$ provides a bijection between
the set of all additive units of $u$ in $(\mathcal E, B)$ and the set of all additive units of $i^*(u)$ in $(E, \beta )$ by letting it acts point-wise on units. More over if $i^*(a)$ is a root of $i^*(u),$ then $a$ is a root of $u.$
\end{prop}
Proof: Suppose $u$ is a unit of the algebraic product system $(\mathcal E,B).$
Then by Theorem \ref{unit}, $i^*(u)$ is a unit of the of
the inclusion system and $\hat{i^*(u)}=u.$ Let $a$ be an additive unit of $u.$
Consider $i^*(a).$ Now \allowdisplaybreaks{\begin{align*} \beta^*_{s,t}[i^*_s(a_s)\otimes
i^*_t(u_t)+i^*_s(u_s)\otimes i^*_t(a_t)] &= [(i_s\otimes
i_t)\beta_{s,t}]^*[a_s\otimes u_t+u_s\otimes a_t]\\ &=
[B_{s,t}i_{s+t}]^*[a_s\otimes u_t+u_s\otimes a_t] \\ &=
i^*_{s+t}a_{s+t}. \end{align*}} %
Hence $i^*(a)$ is an additive unit of the
unit $i^*(u).$

Now we prove the injectivity of $i^*.$ Consider two additive units
$a$ and $b$ of the unit $u$ in $(\mathcal E,B)$ such that
$i^*_ta_t=i^*_tb_t$ for all $t>0.$ Fix $t>0.$ For
$\textbf{s}=(s_1,s_2,...,s_n)\in J_t,$ Define
$a_\textbf{s}=\sum^n_{j=1}u_{s_1}\otimes
u_{s_2}\otimes\cdots u_{s_{j-1}}\otimes a_{s_j}\otimes
u_{s_{j+1}}\otimes\cdots\otimes u_{s_n}$ and
$b_\textbf{s}=\sum^n_{j=1}u_{s_1}\otimes u_{s_2}\otimes \cdots
\otimes u_{s_{j-1}}\otimes b_{s_j}\otimes
u_{s_{j+1}}\otimes\cdots\otimes u_{s_n}.$ Now for $\textbf{s}\in J_t,$

\allowdisplaybreaks{\begin{align*} i^*_\textbf{s}a_t &=
i^*_\textbf{s}B^*_{t,\textbf{s}}a_\textbf{s} \\
&=(B_{t,\textbf{s}}i_\textbf{s})^*a_\textbf{s}\\ &=
(i^*_{s_1}\otimes \cdots \otimes
i^*_{s_n})\sum^n_{j=1} u_{s_1}\otimes u_{s_2}\otimes \cdots
\otimes u_{s_{j-1}}\otimes a_{s_j}\otimes u_{s_{j+1}}\otimes \cdots
\otimes u_{s_n} \\ &= (i^*_{s_1}\otimes \cdots \otimes
i^*_{s_n})\sum^n_{j=1} u_{s_1}\otimes
u_{s_2}\otimes \cdots\otimes u_{s_{j-1}}\otimes b_{s_j}\otimes
u_{s_{j+1}}\otimes \cdots \otimes u_{s_n} \\ &=(B_{t,\textbf{s}}i_\textbf{s})^*b_\textbf{s} \\ &= i^*_\textbf{s}B^*_{t,\textbf{s}}b_\textbf{s} \\ &= i^*_\textbf{s}b_t.\end{align*} } %

This implies $i_\textbf{s}i^*_\textbf{s}a_t=i_\textbf{s}i^*_\textbf{s}b_t.$  The
net of projection $\{i_\textbf{s}i^*_\textbf{s}:\textbf{s}\in J_t\}$
converges strongly to the identity. So we get $a_t=b_t.$

Conversely, let $u$ be a unit and $a$ be an additive unit of $u$ in
$(E,\beta).$ Fix $t>0.$ For $\textbf{s}=(s_1,s_2,...,s_n)\in
J_t,$ Define $a_\textbf{s}=\sum^n_{j=1} u_{s_1}\otimes
u_{s_2}\otimes \cdots \otimes u_{s_{j-1}}\otimes a_{s_j}\otimes
u_{s_{j+1}}\otimes \cdots \otimes u_{s_n}.$ Now the family
$\{i_\textbf{s}a_\textbf{s}:\textbf{s}\in J_t\}$ is bounded as \allowdisplaybreaks{\begin{align*}
\|i_\textbf{s}a_\textbf{s}\|^2 &\leq
\sum^n_{i=1}k(s_i+s_i^2) \\ &\leq k(s+s^2).\end{align*}} %
It follows from the hypothesis that, for $\textbf{s}\leq\textbf{t}\in
J_t,$ \allowdisplaybreaks{\begin{align*} a_\textbf{s} &=
\beta^*_{\textbf{s},\textbf{t}}a_\textbf{t}.\end{align*}} %
Now for $\textbf{s}\leq\textbf{t}\in J_t,$ \allowdisplaybreaks{\begin{align*}
i_\textbf{s}i^*_\textbf{s}i_\textbf{t}a_\textbf{t} =
i_\textbf{s}\beta^*_{\textbf{s},\textbf{t}}a_\textbf{t} =
i_\textbf{s}a_\textbf{s}.\end{align*}} %
Given $\epsilon >0,$ $x\in E_t,$ choose $\textbf{s}\in J_t$ such that $\|(I-i_{\textbf{s}}i^*_{\textbf{s}})x\|<\epsilon.$ Then for $\textbf{t}\geq \textbf{s},$ we have \allowdisplaybreaks{\begin{align*} \langle i_{\textbf{t}}a_{\textbf{t}}-i_{\textbf{s}}a_{\textbf{s}} , x  \rangle  &= \langle (I-i_{\textbf{s}}i^*_{\textbf{s}})i_{\textbf{t}}a_{\textbf{t}},x \rangle \\ &= \langle i_{\textbf{t}}a_{\textbf{t}}, (I-i_{\textbf{s}}i^*_{\textbf{s}})x \rangle \\ & \leq  \|i_{\textbf{t}}a_{\textbf{t}} \| \| (I-i_{\textbf{s}}i^*_{\textbf{s}})x\| \\ & \leq  [k(s+s^2)]^{\frac{1}{2}}\epsilon. \end{align*}} %
So for each $x\in \mathcal E_t,$ $\{\langle i_\textbf{s}a_{\textbf{s}},x\rangle : \textbf{s}\in J_t\}$ is a weakly Cauchy net. Set $\phi(x)=\lim\limits_{\textbf{s}\in J_t} \langle i_\textbf{s}a_{\textbf{s}},x\rangle.$ Then $\phi:\mathcal E_t\rightarrow \mathbb C$ is a bounded linear functional with $\|\phi\|\leq k(s+s^2).$ So there is a unique vector $\hat{a}_t \in \mathcal E_t$ such that $\phi(x)=\langle \hat{a}_t,x\rangle.$ This implies for every $x \in \mathcal E_t,$ $\langle i_\textbf{s}a_{\textbf{s}},x\rangle = \langle \hat{a}_t,x \rangle.$ Now for $\textbf{s}\in J_t,$ \allowdisplaybreaks{\begin{align*} i_{\textbf{s}}i^*_{\textbf{s}}\hat{a}_t &= \lim\limits_{\textbf{t}\in J_t}i_\textbf{s}i^*_\textbf{s}  i_{\textbf{t}}a_{\textbf{t}} \\ &= \lim\limits_{\textbf{t}\in J_t} i_\textbf{s}\beta^*_{\textbf{s},\textbf{t}}a_\textbf{t} \\ &= i_\textbf{s}a_\textbf{s}. \end{align*}} %
This shows that $\{i_{\textbf{s}}a_{\textbf{s}}: \textbf{s} \in J_t\}$ converges to $\hat{a}_t$ in the Hilbert space norm. Let $\hat{u}$ be the lift of $u$ in the algebraic product
system. Our claim is that $\hat{a}=(\hat{a}_t)_{t>0}$ is an additive unit of the unit $\hat{u}=(\hat{u}_t)_{t>0}$
in the algebraic product system. For $x\in\mathcal E_s$, $y\in\mathcal E_t,$
\allowdisplaybreaks{\begin{align*} \langle\hat{a}_s\otimes\hat{u}_t+\hat{u}_s\otimes\hat{a}_t,x\otimes
y\rangle &= \lim\limits_{\textbf{s}\in J_s,\textbf{t}\in J_t}\langle
(i_\textbf{s}\otimes i_\textbf{t})[a_\textbf{s}\otimes
u_\textbf{t}+u_\textbf{s}\otimes a_\textbf{t}],x\otimes y\rangle \\
&= \lim\limits_{\textbf{s}\in J_s,\textbf{t}\in J_t} \langle (i_\textbf{s}\otimes i_\textbf{t}) a_{\textbf{s}\smile\textbf{t}} , (x\otimes y)\rangle \\ &=
\lim\limits_{\textbf{s}\in J_s,\textbf{t}\in J_t}\langle
B_{s,t}i_{\textbf{s}\smile\textbf{t}}a_{\textbf{s}\smile\textbf{t}},x\otimes
y\rangle \\ &= \langle \lim\limits_{\textbf{s}\smile\textbf{t} \in J_{s+t}} i_{\textbf{s}\smile\textbf{t}} a_{\textbf{s}\smile\textbf{t}} , B^*_{s,t}(x\otimes y) \rangle \\ &= \langle B_{s,t} \hat{a}_{s+t},x\otimes y \rangle. \end{align*}} %
This proves the claim.

For $x\in E_t,$ we have \allowdisplaybreaks{\begin{align*} \langle
i^*_t\hat{a}_t,x\rangle &= \langle\hat{a}_t,i_tx\rangle \\ &=
\lim\limits_{\textbf{r}\in J_t}\langle i_\textbf{r}a_\textbf{r},i_tx\rangle \\
&= \lim\limits_{\textbf{r}\in J_t}\langle i^*_ti_\textbf{r}a_\textbf{r},x\rangle \\ &=
\lim\limits_{\textbf{r}\in J_t}\langle
\beta^*_{t,\textbf{r}}i^*_\textbf{r}i_\textbf{r}a_\textbf{r},x\rangle
\\ &= \lim\limits_{\textbf{r}\in J_t}\langle \beta^*_{t,\textbf{r}}a_\textbf{r},x\rangle
\\ &= \langle a_t,x\rangle. \end{align*}} %
This implies $i^*_t\hat{a}_t=a_t.$

Finally, if $b$ is a root of a unit $v$ in the inclusion system $(E,\beta),$ then \allowdisplaybreaks{\begin{align*} \langle \hat{b}_t,\hat{v}_t \rangle & = \lim\limits_{\textbf{r}\in J_t} \langle i_{\textbf{r}}b_{\textbf{r}}, i_{\textbf{r}}v_{\textbf{r}} \rangle \\ &= \lim\limits_{\textbf{r}\in J_t} \langle b_{\textbf{r}} , v_{\textbf{r}} \rangle \\ &= \lim\limits_{\textbf{r}\in J_t} \sum^n_{j} \langle v_{r_1}\otimes
v_{r_2}\otimes\cdots v_{r_{j-1}}\otimes b_{r_j}\otimes
v_{r_{j+1}}\otimes\cdots\otimes v_{r_n} , v_{r_1}\otimes
v_{r_2}\otimes\cdots\otimes v_{r_n}\rangle \\ & =  0.
 \end{align*}} %

This proves the last assertion.
\qed

Here we show how the root space behaves under the amalgamation
through partial isometry. Suppose $\mathcal E$ and $\mathcal F$ are two product systems and
$C=(C_t)_{t>0}:\mathcal F_t\rightarrow \mathcal E_t$ is a morphism
of partial isometry. Also assume that $\mathcal E\otimes_C\mathcal F$ is a product system. Let $v=(v_t)_{t>0}$ be a normalized unit of $\mathcal F$ such that $C^*_tC_tv_t=v_t$ for all $t>0.$ Then $Cv:=(C_tv_t)_{t>0}$ is a normalized unit of $\mathcal E.$  We have for every $t>0,$ $I_t{C_tv_t}=J_t{v_t},$ where $I_t:\mathcal E_t\rightarrow (\mathcal E\otimes_C\mathcal F)_t$ and $J_t:\mathcal F_t\rightarrow (\mathcal E\otimes_C\mathcal F)_t$ are injection morphisms. Let us denote the common unit by $u$ in $\mathcal E\otimes_C\mathcal F.$ i.e. $u_t=I_t(C_tv_t)=J_t(v_t)$ for all $t>0.$
Denoting $\mathcal R^\mathcal E_{Cv}:= \{a_1\in\mathcal E_1:a\in R^\mathcal E_{Cv}\}$ and $\mathcal R^\mathcal F_v:= \{b_1\in\mathcal F_1:b\in R^\mathcal F_v\}.$ Given two closed subspaces $H$ and $H^\prime$ of a Hilbert space $G,$ denote by $H\vee H^\prime$ the smallest closed subspace of $G$ containing $H$ and $H^\prime.$

\begin{thm} \label{partialroot}
 Suppose ($\mathcal E,W^{\mathcal E})$ and
 ($\mathcal F,W^{\mathcal F}$) are two spatial product systems and suppose
 $C=(C_t)_{t>0}:\mathcal F_t\rightarrow \mathcal E_t$ is a morphism
 of partial isometry such that the amalgamated product $\mathcal E\otimes_C \mathcal F$ is a product system. Also Suppose $v=(v_t)_{t>0}$ is a normalized unit of $\mathcal F$ such that $C^*_tC_tv_t=v_t$ for all $t>0.$ Then $Cv:=(C_tv_t)_{t>0}$ and $v=(v_t)_{t>0}$ are identified in $\mathcal E\otimes_C\mathcal F.$ Denote the common unit by $u=(u_t)_{t>0}$ in $\mathcal E\otimes \mathcal F.$ Then $\mathcal
R^{\mathcal E\otimes_C\mathcal F}_u=\mathcal R^{\mathcal E}_{Cv}
\oplus_{C_1} \mathcal R^{\mathcal F}_v.$
\end{thm}

Proof:  We may assume from Theorem 2.7, \cite{Muk-index}, that $\mathcal E$ and $\mathcal F$ are subsystems
of the amalgamated product $\mathcal E\otimes_C\mathcal F.$ As $C$
is a morphism of partial isometry, we get from \cite{Muk-index}, Proposition 2.10, that for each $t>0,$
$P_{\mathcal E_t}$ and $P_{\mathcal F_t}$ commute as elements in $\mathcal B((\mathcal E\otimes_C\mathcal F)_t).$ So $P_{\mathcal E_t\cap
\mathcal F_t}=P_{\mathcal E_t}P_{\mathcal F_t},$ which implies  $\mathcal E\cap\mathcal F:= (\mathcal E_t\cap\mathcal F_t)_{t>0}$ is a product subsystem. In this identification, we have $u=v=Cv.$ Hence $u$ is a normalized unit of $\mathcal E\cap \mathcal F$ and $\mathcal R^{\mathcal E}_u\oplus_{C_1} R^{\mathcal F}_u$ coincides with $\mathcal R^{\mathcal E}_u\vee \mathcal R^{\mathcal F}_u$ inside $\mathcal R^{\mathcal E\otimes_C\mathcal F}_u.$ So to prove the theorem, it is enough to show that $$\mathcal R^{\mathcal E}_u\vee \mathcal R^{\mathcal F}_u=\mathcal R^{\mathcal E\otimes_C \mathcal F}_u.$$ Clearly $\mathcal R^{\mathcal E}_u\vee \mathcal R^{\mathcal F}_u \subset \mathcal R^{\mathcal E\otimes_C \mathcal F}_u.$ Now for $a\in R^{\mathcal E\otimes_C \mathcal F}_u,$ consider $b=(b_t)_{t>0}$ where $b_t=P_{\mathcal E_t}a_t,$  $b^\prime=(b^\prime_t)_{t>0}$ where $b^\prime_t=P_{\mathcal F_t}a_t$ and $b^{\prime\prime}=(b^{\prime\prime}_t)_{t>0}$ where $b^{\prime\prime}_t=P_{\mathcal E_t\cap \mathcal F_t}a_t.$
We claim
that $$b\in R^{\mathcal E}_u,~ b^\prime \in R^{\mathcal F}_u,~ b^{\prime\prime}\in R^{\mathcal E\cap \mathcal F}_u.$$  Note that for every $s>0,$ $$u_s=P_{\mathcal E_s}=P_{\mathcal F_s}u_s=P_{\mathcal E_s\cap\mathcal F_s}u_s.$$ As $P_{\mathcal E}=(P_{\mathcal E_s})_{s>0}$ is a projection morphism from $(\mathcal E\otimes_C\mathcal F, W^{\mathcal E\otimes \mathcal F})$ to $(\mathcal E,W^{\mathcal E}),$ we have $$ (P_{\mathcal E_s}\otimes P_{\mathcal E_t})W^{(\mathcal E\otimes_C\mathcal F)}_{s+t} =  W^{\mathcal E}_{s,t} P_{\mathcal E_{s+t}},~s,t>0. $$ This implies \allowdisplaybreaks{\begin{align*} W^{\mathcal E}_{s,t} b_{s+t} &= W^{\mathcal E}_{s,t} P_{\mathcal E_{s+t}} a_{s+t} \\ &= (P_{\mathcal E_s}\otimes P_{\mathcal E_t})W^{(\mathcal E\otimes_C\mathcal F)}_{s+t} a_{s+t} \\ &= (P_{\mathcal E_s}\otimes P_{\mathcal E_t})(a_s\otimes u_t+u_s\otimes a_t) \\ &= (b_s\otimes u_t + u_s \otimes b_t). \end{align*}} %
This shows $b\in R^{\mathcal E}_u.$ Similarly we have, $b^\prime\in R^{\mathcal F}_u$ and $b^{\prime\prime}\in R^{\mathcal E\cap\mathcal F}_u.$ Also note that $b,b^\prime,b^{\prime\prime} \in R^{\mathcal E\otimes_C\mathcal F}_u.$
Set $c=b+b^\prime-b^{\prime\prime}.$ Then we have for all $t>0,$ $$P_{\mathcal E_t}(a_t-c_t)=b_t-b_t=0.~P_{\mathcal F_t}(a_t-c_t)=0,.$$ Therefore $$P_{\mathcal E_t\vee\mathcal F_t}(a_t-c_t)=0.$$ Note that $(\mathcal E_t\vee\mathcal F_t)_{t>0})$ is an inclusion system which generates the product system $\mathcal E\otimes_C\mathcal F.$ Also note that $(P_{\mathcal E_t\vee\mathcal F_t}(a_t-c_t))_{t>0}$ is a root of $u$ in the inclusion system $(\mathcal E_t\vee\mathcal F_t)_{t>0})$ while $(a_t-c_t)_{t>0}$ is a root of $u$ in the product system $(\mathcal E \otimes_C\mathcal F).$ As $(\mathcal E_t\vee\mathcal F_t)_{t>0})$ generates the product system $(\mathcal E \otimes_C\mathcal F),$ we have from the injectivity of the map $i^*$ described in Theorem \ref{root},  for all $t>0,$ $a_t=c_t.$  So $a_1=b_1-b^{\prime\prime}_1+b^\prime_1,$ where $b_1-b^{\prime\prime}_1\in \mathcal R^{\mathcal E}_u$ and $b^\prime_1 \in \mathcal R^{\mathcal F}_u.$ Hence $\mathcal R^{\mathcal E}_u\vee \mathcal R^{\mathcal F}_u=\mathcal R^{\mathcal E\otimes_C \mathcal F}_u.$ \qed

Suppose ($\mathcal E,W^{\mathcal E})$ and ($\mathcal F,W^{\mathcal
F}$) are two product systems. Let $u^0$ and $v^0$ be two normalized
units of $\mathcal E$ and $\mathcal F$ respectively. Consider
$\mathcal E\otimes_C \mathcal F,$ where $C_t=|u^0_t\rangle\langle
v^0_t|.$ In the amalgamated product system $\mathcal
E\otimes_C\mathcal F,$ $u^0$ and $v^0$ are identified. We denote the
common unit by $\sigma.$

\begin{cor}
Let $\mathcal E,$ $\mathcal F,$ $u^0,$ $v^0,$ $\sigma$ be as above.
Then $\mathcal R^{\mathcal E\otimes_C\mathcal F}_\sigma=\mathcal
R^{\mathcal E}_{u_0}\oplus \mathcal R^{\mathcal F}_{v_0}.$
\end{cor}

Proof: For $x\in \mathcal R^{\mathcal
E}_{u_0},$ $y\in \mathcal R^{\mathcal F}_{v_0},$ \allowdisplaybreaks{\begin{align*} \langle x,y
\rangle_{C_1} &= \langle x,C_1y \rangle \\ &= \langle x,
|u^0_1\rangle \langle v^0_1|y \rangle \\ &= \langle x,u^0_1
\rangle \langle v^0_1,y \rangle \\ &= 0 . \end{align*}} %
\qed

\begin{rem}
It is noted that the condition on $C$ that it is a partial isometry
in Theorem \ref{partialroot} is a necessary condition. It may not be
true for general contractive morphism. Let $\mathcal E_t=\mathbb
Cu_t$ and $\mathcal F_t=\mathbb Cv_t$ be two type $I_0$ product systems
with $\|u_t\|\|v_t\| < 1$  for some $t>0.$ Let
$C_t=|u_t\rangle\langle v_t|.$ Then $\mathcal R^{\mathcal E}_u=0$
and $\mathcal R^{\mathcal F}_v=0.$ On the other side, we have
$\mathcal E\otimes_C\mathcal F$ is a type $I_1$ product system. Though a priori, it is not clear whether in this case, $\mathcal E\otimes_C\mathcal F$ is a product system. But this is indeed true (\cite{Muk-contractive}).
Therefore $\mathcal R^{\mathcal E\otimes_C\mathcal F}_\sigma \neq
\{0\}$ for every unit $\sigma$ in $\mathcal E\otimes_C\mathcal F.$
Hence $\mathcal R^{\mathcal E\otimes_C\mathcal F}_u \neq \mathcal
R^{\mathcal E}_u \oplus_{C_1} \mathcal R^{\mathcal F}_u.$
\end{rem}

\section{Amalgamation through normalized units is independent of the choice of units }

In this section, we will show that the amalgamation through normalized
unit does not depend on the choice of the units. Proof of this fact is
almost visible when we use the theory of random sets (\cite{Lie-random}). In \cite{BLMS-intrinsic}, a short and self-contained proof has been presented. Also see \cite{BLS-Powers}. Here we will prove this fact using roots.

First, we show that the amalgamation of two spatial product systems
through normalized units can be identified with the product
subsystem of the tensor product of the two systems. Let $\mathcal E$
and $\mathcal F$ be two spatial product systems and $u$ and $v$ be two normalized
units of $\mathcal E$ and $\mathcal F$ respectively.
Define a contractive morphism $C=(C_t)_{t>0}:\mathcal F_t\rightarrow
\mathcal E_t$ by $C_t=|u_t\rangle\langle v_t|.$ Denote $\mathcal
E\otimes_{u,v}\mathcal F:=\mathcal E\otimes_C\mathcal F.$ For two product subsystems $\mathcal G$ and $\mathcal G^\prime$ of the product system $\mathcal H,$ we denote by $\mathcal G\bigvee\mathcal G^\prime$ the smallest product subsystem of $\mathcal H$ containing $\mathcal G$ and $\mathcal G^\prime.$

\begin{prop}\label{tensor}
Suppose $\mathcal E$ and $\mathcal F$ are two spatial product systems and
$u$ and $v$ are two normalized units of $\mathcal E$ and $\mathcal F$
respectively. Then $\mathcal E\otimes_{u,v}\mathcal F$ is isomorphic
to the product system generated by $\mathcal E\otimes v$ and
$u\otimes\mathcal F$ inside $\mathcal E\otimes \mathcal F,$ i.e.
$\mathcal E\otimes_{u,v}\mathcal F\simeq  (\mathcal E\otimes
v)\bigvee (u\otimes \mathcal F).$
\end{prop}
Proof: As $u$ and $v$ are normalized, we see that $I:\mathcal E\rightarrow \mathcal E\otimes v $ and $J:\mathcal F\rightarrow u\otimes\mathcal F$ are isometric morphisms of product system. Also note that for $x \in \mathcal E_s$ and $y\in \mathcal F_s,$  $\langle I(x),J(y)\rangle=\langle x, |u_t\rangle\langle v_t|y \rangle.$ Now from the property of amalgamation (Theorem 2.7, \cite{Muk-index}) we conclude that $\mathcal E\otimes_{u,v}\mathcal F\simeq  (\mathcal E\otimes v) \bigvee (u\otimes
                       \mathcal F) \subset \mathcal E\otimes \mathcal F$ as algebraic product systems. Now transferring the measurable structure of $ (\mathcal E\otimes v) \bigvee (u\otimes \mathcal F)$ onto $\mathcal E\otimes_{u,v}\mathcal F $ via the isomorphism, we can make $\mathcal E\otimes_{u,v}\mathcal F$ into a product system and the isomorphism becomes the isomorphism of product systems.  \qed

Suppose $\mathcal E$ is a product system and $u=(u_t)_{t>0}$ is a normalized unit of $\mathcal E.$  Then for every interval $[s,t],$ $0<s<t<1,$ we may identify, $\mathcal E_1\simeq \mathcal E_{s}\otimes\mathcal E_{t-s}\otimes\mathcal E_{1-t}.$   Let $P_{s,t}= P_{\mathcal E_{s}\otimes\mathbb Cu_{t-s}\otimes\mathcal E_{1-t}} = 1_{\mathcal E_s}\otimes P_{\mathbb C u_{t-s}}\otimes
1_{\mathcal E_{1-t}}.$ From Proposition 3.18, \cite{Lie-random}, we know that $(s,t)
\rightarrow P_{s,t}$ is jointly continuous. So in the compact
simplex $\{0\leq s\leq t \leq 1\},$ it is uniformly continuous. i.e.
$P_{s,t}$ goes to identity strongly as $(t-s) \rightarrow 0.$ In this section, we denote the multiplication operation of the product system by $\circ$ i.e. $a\in \mathcal E_s,$ $b\in\mathcal E_t,$ we have $a\circ b\in \mathcal E_{s+t}.$ We write $P_{s,t}$ as  $1_{\mathcal E_s}\circ P_{\mathbb C u_{t-s}}\circ
1_{\mathcal E_{1-t}}.$ This is to differentiate the multiplication operation of the product system with the tensor product operation on the category of product systems. Though note that this is not the usual operator multiplications as they are not acting on the same space. We hope these notations do not lead any confusion.

For $n\geq 1,$ we have $P_{\frac{i-1}{n},\frac{i}{n}}=1_{\mathcal E_{\frac{1}{n}}}\circ
\cdots\circ 1_{\mathcal E_{\frac{1}{n}}}\circ P_{\mathbb C u_{\frac{1}{n}}} \circ
\cdots \circ 1_{\mathcal E_{\frac{1}{n}}},$ where $P_{\mathbb C u_{\frac{1}{n}}}$
on the $i$-th place.
\begin{thm}
Suppose $\mathcal E$ and $\mathcal F$ are two spatial product
systems with normalized units $u$ and $v$ respectively. Then
$\mathcal E\otimes_{u,v}\mathcal F$ is isomorphic to the product
system generated by $\mathcal E\otimes \mathcal F^I$ and $\mathcal
E^I\otimes \mathcal F$ inside $\mathcal E\otimes \mathcal F,$ i.e.
$\mathcal E\otimes_{u,v}\mathcal F\simeq  (\mathcal E\otimes
\mathcal F^I) \bigvee (\mathcal E^I\otimes\mathcal F).$
\end{thm}
Proof: We know from Proposition \ref{tensor}, that $\mathcal
E\otimes_{u,v}\mathcal F\simeq  (\mathcal E\otimes v)\bigvee
(u\otimes \mathcal F) \subset \mathcal E \otimes \mathcal F.$ So to
prove the theorem, it is enough to show that $\mathcal E\otimes
\mathcal F^I\subset (\mathcal E\otimes v)\bigvee (u\otimes \mathcal
F),$ as the proof of $\mathcal E^I\otimes \mathcal F\subset
(\mathcal E\otimes v)\bigvee (u\otimes \mathcal F)$ is identical. We
fix the time point $t=1.$ Now from Theorem \ref{unit-roots}, it is
enough to show that for $z\in \mathcal E_1$ and  for any root $a$ of
$v$ with $\|a_1\|=1,$ $z\otimes a_1\in ((\mathcal E\otimes v)\bigvee
(u\otimes \mathcal F))_1.$ For other time point, proof goes
identically. Let $\epsilon > 0$ be given. From uniform continuity
of $P_{s,t},$ choose $N$ such that $n \geq N,$ $\|
z-P_{\frac{i-1}{n},\frac{i}{n}}z\| \leq \epsilon,$ for every $i=1,2,\cdots,n.$
Choose and fix $n\geq N.$ Decompose $a_1=\sum_{i=1}^n x^i,$ where
$x^i=v_{\frac{1}{n}}\circ\cdots\circ v_{\frac{1}{n}}\circ a_{\frac{1}{n}}\circ
v_{\frac{1}{n}}\circ\cdots\circ v_{\frac{1}{n}},$ with $a_{\frac{1}{n}}$ at i-th place.
Clearly $\|x^i\|=1/\sqrt{n}.$ \allowdisplaybreaks{\begin{align*} \|z\otimes a_1-\sum_{i=1}^n
(P_{\frac{i-1}{n},\frac{i}{n}}z\otimes x^i)\|^2 &=
\|\sum_i^n z\otimes x^i-\sum_i^n (P_{\frac{i-1}{n},\frac{i}{n}} z\otimes x^i)\|^2 \\
&= \|\sum_i^n(z-P_{\frac{i-1}{n},\frac{i}{n}}z)\otimes x_i\|^2 . \end{align*}} %
Now for $i \neq j$ we have $\langle x_i,x_j \rangle =0.$ So $\langle
(z-P_{\frac{i-1}{n},\frac{i}{n}}z)\otimes x_i,(z-P_{\frac{j-1}{n},\frac{j}{n}}z)\otimes x_j
\rangle=0.$ So we have \allowdisplaybreaks{\begin{align*} \|z\otimes a_1-\sum_{i=1}^n
(P_{\frac{i-1}{n},\frac{i}{n}}z\otimes x^i)\|^2 &=
\sum_i^n\|z-P_{\frac{i-1}{n},\frac{i}{n}}z \|^2\|x^i\|^2 \\ &\leq 1/n \sum_i^n\|z-P_{\frac{i-1}{n},\frac{i}{n}}z \|^2 \\ &\leq \epsilon^2. \\ \end{align*}} %
Now the vector \allowdisplaybreaks{\begin{align*}  P_{\frac{i-1}{n},\frac{i}{n}} z\otimes x^i &=
\sum_j(c^1_j\circ c^2_j\circ\cdots\circ c^{i-1}_j\circ u_{\frac{1}{n}}\circ
c^{i+1}_j\circ\cdots\circ c^n_j )\\ &\otimes(v_{\frac{1}{n}}\circ
v_{\frac{1}{n}}\circ\cdots\circ a_{\frac{1}{n}}\circ v_{\frac{1}{n}}\circ\cdots\circ
v_{\frac{1}{n}})\\ &= \sum_j (c^1\otimes v_{\frac{1}{n}})\circ (c^2\otimes
v_{\frac{1}{n}})\circ\cdots\circ (c^{j-1}\otimes v_{\frac{1}{n}})\\ & \circ
(u_{\frac{1}{n}}\otimes a_{\frac{1}{n}}) \circ
(c^{j+1}\otimes v_{\frac{1}{n}})\circ\cdots\circ (c^n\otimes v_{\frac{1}{n}}).\\ \end{align*}} %

Note that, for every $1\leq j\leq n,$ $c^j\otimes v_{\frac{1}{n}} \in (\mathcal E \otimes v)_{\frac{1}{n}}$ and $(u_\frac{1}{n}\otimes a_\frac{1}{n})\in (u\otimes\mathcal F)_\frac{1}{n}.$

This implies that $z\otimes a_1 \in ((\mathcal E\otimes v)\bigvee
(u\otimes \mathcal F))_1.$ \qed

\begin{cor}\label{root-typeI}
Suppose $\mathcal E$ and $\mathcal F$ are two spatial product systems  with normalized units $u$ and $v$ respectively. Then $(\mathcal E^I\otimes v)\bigvee(u\otimes \mathcal F^I)=(\mathcal E^I\otimes\mathcal F^I)=(\mathcal E\otimes \mathcal F)^I.$
\end{cor}

\begin{thm}
Suppose $\mathcal E$ and $\mathcal F$ are two spatial product systems  with normalized units $u$ and $v$ respectively. Then $R^{\mathcal E\otimes\mathcal F}_{u\otimes v}=(R^{\mathcal E}_u\otimes v) \oplus (u\otimes R^{\mathcal F}_v).$
\end{thm}
Proof: For $a\in R^{\mathcal E}_u$ and $b\in R^{\mathcal F}_v,$ define for each $s>0,$ $d_s=a_s\otimes v_s + u_s\otimes b_s.$ Then for $s,t>0,$ we have   \allowdisplaybreaks{\begin{align*} & d_s\circ(u_t\otimes v_t)+(u_s\otimes v_s)\circ d_t \\ &= (a_s\otimes v_s + u_s\otimes b_s)\circ(u_t\otimes v_t) + (u_s\otimes v_s)\circ(a_t\otimes v_t + u_t\otimes b_t) \\ &= [(a_s\circ u_t)\otimes(v_s\circ v_t)+(u_s\circ u_t)\otimes(b_s\circ v_t)]+ [(u_s\circ a_t)\otimes(v_s\circ v_t)+ (u_s\circ u_t)\otimes (v_s\circ b_t)] \\ &= [(a_s\circ u_t+u_s\circ a_t)\otimes v_{s+t}+u_{s+t}\otimes (b_s\circ v_t+v_s\circ b_t)] \\ &= a_{s+t}\otimes v_{s+t}+u_{s+t}\otimes b_{s+t} \\ &= d_{s+t}. \end{align*} } %
This implies $d \in R^{\mathcal E\otimes \mathcal F}_{u\otimes v}$ and we obtain $R^{\mathcal E\otimes\mathcal F}_{u\otimes v} \supset (R^{\mathcal E}_u\otimes v) \oplus (u\otimes R^{\mathcal F}_v).$ Observe that $\{(a_s\otimes v_s) + (u_s\otimes b_s):a\in R^{\mathcal E}_u,b\in R^{\mathcal F}_v\}$ is a closed subspace of $\{c_s: c\in R^{\mathcal E\otimes\mathcal F}_{u\otimes v}\}.$ We  claim that for every $s>0,$  $\{c_s: c\in R^{\mathcal E\otimes\mathcal F}_{u\otimes v}\} \ominus \{(a_s\otimes v_s) + (u_s\otimes b_s):a\in R^{\mathcal E}_u,b\in R^{\mathcal F}_v\}=0.$ Suppose $c\in R^{\mathcal E\otimes \mathcal F}_{u\otimes v}.$ Also suppose that for all $a\in R^{\mathcal E}_u$ and $b\in R^{\mathcal F}_v,$ and for every $s>0,$  $\langle c_s,a_s\otimes v_s+u_s\otimes b_s\rangle=0.$  As $c\in R^{\mathcal E\otimes \mathcal F}_{u\otimes v},$ we have for every $s>0,$ $\langle c_s, u_s\otimes v_s\rangle=0.$ Now from Corollary \ref{unit-roots}, we have for every $s>0,$ $c_s$ belong to the ortho-complement of $((\mathcal E^I\otimes v) \bigvee (u\otimes \mathcal F^I))_s.$ Now from Corollary \ref{root-typeI}, we get, for every $s>0,$ $c_s$ is in the ortho-complement of $(\mathcal E\otimes \mathcal F)^I_s.$ But as $c\in R^{\mathcal E\otimes \mathcal F}_{u\otimes v},$ we have from Corollary \ref{root I}, for every $s>0,$ $c_s \in (\mathcal E\otimes \mathcal F)^I_s.$ This shows that for every $s>0,$ $c_s=0.$ This proves the claim. Hence we have the equality $R^{\mathcal E\otimes\mathcal F}_{u\otimes v}=(R^{\mathcal E}_u\otimes v) \oplus (u\otimes R^{\mathcal F}_v).$     \qed

\section{Cluster construction}
Here we introduce a new construction called the cluster construction.
Given any product subsystem $\mathcal F$ in a product system $\mathcal E,$
we attach a product subsystem $\check{\mathcal F}\supset \mathcal F.$ We call the
product subsystem $\check{\mathcal F}$ as the cluster of $\mathcal F$ in $\mathcal E.$
The name `cluster' comes from the following connection of random
sets discussed in \cite{Lie-random}. Every product subsystem corresponds to a unique probability measure on the closed subsets of $[0,1].$ The set of all closed sets of $[0,1]$ can be topologized by hit and miss topology (see Page 2, \cite{Lie-random} , Section 1-4, \cite{Mat-random} for details). The mapping `cluster' which sends a closed set to its limit points is a measurable map on this space. We show here that the probability measure corresponding to the cluster subsystem is the distribution of the cluster map. We compute the `cluster' of the product subsystem of a spatial product system given by a single unit and show that it is the type I part of the product system.

Suppose $(\mathcal E,B)$ is a product system and $(F,B|_F)$ is an
inclusion subsystem. Define $\tilde{F_t}$ by
$$\tilde{ F_t}=\overline{\mbox{span}}\{x\otimes y: x \in \mathcal E_r\ominus F_r, y \in \mathcal E_{t-r}\ominus F_{t-r}, ~\mbox{for}~ \mbox{some}~ r,
0<r<t\}. $$ Set $F^\prime_t=\mathcal E_t\ominus \tilde{F_t}.$
\begin{lem}
With the notation as above, $(F^\prime_t,B_{s,t}|_{F^\prime_t})$ is
an inclusion system.
\end{lem}
Proof: Let $x \in F^\prime_{s+t}.$ First note that $$F^\prime_s\otimes F^\prime_t=(\mathcal E_s\otimes F^\prime_t)\cap (F^\prime_s \otimes \mathcal E_t).$$ Now $$F_{s+r} \subset F_s\otimes F_r \subset \mathcal E_s\otimes F_r,$$ implies $$\mathcal E_s\otimes (\mathcal E_r\ominus F_r)\subset \mathcal E_{s+r}\ominus F_{s+r}.$$ Now For $y \in \mathcal E_s,$ $z_1\in \mathcal E_r\ominus F_r,$ $z_2\in \mathcal E_{t-r}\ominus F_{t-r},$ for some $0<r<t,$ we get $y\otimes z_1 \in  \mathcal E_{s+r}\ominus F_{s+r}.$ So $\langle x, y_1\otimes z_1\otimes z_2 \rangle =0.$ This shows $$x \in \mathcal E_{s+t}\ominus (\mathcal E_s\otimes \tilde{F_t}).$$ i.e. $x\in \mathcal E_s\otimes F^\prime_t.$ Similarly we get for $0 <r^\prime < s,$ $$(\mathcal E_{s-r^\prime}\ominus F_{s-r^\prime})\otimes \mathcal E_t \subset \mathcal E_{s+t-r^\prime}\ominus F_{s+t-r^\prime}.$$ So for $z^\prime_1 \in \mathcal E_{r^\prime}\ominus F_{r^\prime},$ $z^\prime_2\in \mathcal E_{s-r^\prime}\ominus F_{s-r^\prime},$ $y^\prime \in \mathcal E_t,$ we have $z^\prime_2\otimes y^\prime \in \mathcal E_{s+t-r^\prime}\ominus F_{s+t-r^\prime}.$ This shows $$x\in \mathcal E_{s+t}\ominus (\tilde{F_s}\otimes \mathcal E_t).$$ i.e. $x \in F^\prime_s\otimes \mathcal E_t.$ Associativity property follows from the associativity of the product system. \qed

Given a product subsystem $\mathcal F$ of a product system $\mathcal E,$ denote by $\check{\mathcal F}$ the product system generated by the inclusion
system $(\mathcal F^\prime,W|_{\mathcal F^\prime}).$ We call this product subsystem as the cluster of $\mathcal F$ in $\mathcal E.$ Now our present task is to relate the cluster construction with the theory of random sets described in \cite{Lie-random}. Recall that the random closed sets are characterized by the random variables $X_{s,t}=\chi_{\{Z:Z\cap[s,t]=\emptyset\}}(Z),$ fulfilling $X_{r,s}X_{s,t}=X_{r,t},$ $0\leq r \leq s \leq t \leq 1.$  Theorem 3.16, [Lie] shows that the embedding of the product subsystem into the whole product system, i.e. the structure encoded in the algebraic properties of the projections $P^{\mathcal F}_{s,t}$ uniquely determines a measure type of the random closed sets of $[0,1].$ These results translate some of the structure theory of product systems to the structure theory of measure types on the closed subsets of $[0,1].$ Suppose $\mathcal
F$ is a product subsystem of the product system $\mathcal E.$ Fix a faithful normal state $\eta$ on
$B(\mathcal E_1).$ Suppose $\mu^{\mathcal F}_{\eta}$ is the unique probability measure on $\mathfrak F_I,$ as in Theorem 3.16, \cite{Lie-random}, i.e. \bea \label{mu} \mu^{\mathcal F}_\eta\{Z:Z\cap[s_i,t_i]=\emptyset,
i=1,2,\cdots,k\}=\eta(P^{\mathcal F}_{s_1,t_1}\cdots P^{\mathcal
F}_{s_k,t_k}),{~}{~}((s_i,t_i)\in I).\eea  Moreover the correspondence \bea
\label{injection}\chi_{\{Z:Z\cap [s,t]=\emptyset\}}\rightarrow
P^{\mathcal F}_{s,t},{~}((s,t)\in I) \eea extends to an injective
normal representation $J^{\mathcal F}_{\eta}$ of
$L^\infty(\mu^{\mathcal F}_\eta)$ on $\mathcal E_1.$ Its image is
$\{P^{\mathcal F}_{s,t}(s,t)\in I\}^{\prime\prime}.$
For any $Z\in \mathfrak F_K,$ denote $\check{Z}$ the set of its cluster points: $$\check{Z}=\{t\in Z:t \in \overline{ Z\setminus\{t\}}\}.$$ Suppose $l:\mathfrak F_K\rightarrow \mathfrak F_K$ is the measurable map defined by $l(Z)=\check{Z}.$ With these preparations in hand, we can derive an interesting relation between measure types $\mathcal M^{\mathcal F}$ and $\mathcal M^{\check{\mathcal F}}.$
\begin{thm}
Suppose $\mathcal E,$ $\mathcal F,$ $\mu^{\mathcal F}_{\eta}$ are as
above. Then \bea \label{check} J^{\mathcal F}_\eta(\chi_{\{Z:\check{Z}\cap [s,t]=\emptyset\}})=P^{\check{\mathcal F}}_{s,t},{~}((s,t)\in I).\eea Therefore \bea \label{hat} \mu^{\mathcal
F}_{\eta}\{Z:\check{Z}\cap [s_i,t_i]=\emptyset,
i=1,2,\cdots,k\}=\eta(P^{\check{\mathcal F}}_{s_1,t_1}\cdots
P^{\check{\mathcal F}}_{s_k,t_k}) {~}{~}((s_i,t_i)\in I).\eea
Consequently, $\mathcal M^{\check{\mathcal F}}=\mathcal M^{\mathcal F}\circ l^{-1}.$
\end{thm}
Proof: First note that, Equation (\ref{check}) implies that $P^{\check{\mathcal F}} \in \{P^{\mathcal F}_{s,t}:(s,t)\in I\}^{\prime\prime}.$ Now it is enough to prove Equation (\ref{check}).
Indeed \allowdisplaybreaks{\begin{align*} J^{\mathcal F}_{\eta}(\chi_{\{Z:\check{Z}\cap
[s_i,t_i]=\emptyset, i=1,2,\cdots,k\}}) &= J^{\mathcal
F}_{\eta}(\Pi^k_{i=1}(\chi_{\{Z:\check{Z}\cap
[s_i,t_i]=\emptyset\}})
\\ &= \Pi^k_{i=1}J^{\mathcal
F}_{\eta} (\chi_{\{Z:\check{Z}\cap [s_i,t_i]=\emptyset\}}) \\ &=
\Pi^k_{i=1}P^{\check{\mathcal F}}_{s_i,t_i}.\end{align*}} %
Now applying the states $\int (\cdot) d\mu^{\mathcal F}_\eta$ and $\eta$ on $L^\infty(\mu^{\mathcal F}_\eta)$ and $\{P^{\mathcal F}_{s,t}:(s,t)\in I\}^{\prime\prime}$ respectively, we obtain Equation (\ref{hat}).
Now for a closed random set $Z,$ and the interval $[s,t],$ if we have
 $\sharp\{Z\cap[s,t]\}\geq 2.$  i.e. if $Z$ intersects $[s,t]$ at more than one point, then there
exists a rational $q\in \mathbb Q,$ such that $Z \cap [s,q)\neq
\emptyset$ and $Z\cap [q,t]\neq \emptyset.$ So we have the identity,\allowdisplaybreaks{\begin{align*} \{Z:\sharp\{Z\cap [s,t]\}\leq1\} &= [\cup_{q\in \mathbb
Q\cap(s,t)}\{Z:Z\cap [s,q)\neq\emptyset,Z\cap
[q,t]\neq\emptyset\}]^c.
 \end{align*}} %
 Applying $J^{\mathcal F}_{\eta}$ on the indicator function of the above two sets, we get \allowdisplaybreaks{\begin{align*} J^{\mathcal F}_{\eta}(\chi_{ \{Z:\sharp\{Z\cap [s,t]\}\leq1\}}) &= J^{\mathcal F}_{\eta} (\chi_{[\cup_{q\in
\mathbb Q\cap(s,t)}\{Z:Z\cap [s,q)\neq\emptyset,Z\cap
[q,t]\neq\emptyset\}]^c}). \end{align*}} %
Now let $A_q=\{Z:Z\cap [s,q)\neq \emptyset ] \},$ $B_q=\{Z:Z\cap [q,t]\neq \emptyset\}.$ Then
${A^c_q}=\cap_n\{Z:Z\cap [s,q-1/n]=\emptyset\}.$ Further from Equation
\ref{injection} and  continuity of $P^{\mathcal F}_{s,t},$ Proposition 3.18, \cite{Lie-random},  we get $J^{\mathcal
F}_{\eta}(\chi_{{A^c_q}})=\mbox{lim}_n P^{\mathcal
F}_{s,q-1/n}=P^{\mathcal F}_{s,q}.$ Hence $$J^{\mathcal
F}_{\eta}(\chi_{A_q})=1_{\mathcal E_s}\otimes P_{\mathcal
F^\bot_{q-s}}\otimes 1_{\mathcal E_{1-q}}.$$ Similarly,
$$J^{\mathcal F}_{\eta}(\chi_{B_q})=1_{\mathcal E_q}\otimes
P_{\mathcal F^\bot_{t-q}}\otimes 1_{\mathcal E_{1-t}}.$$ Now \allowdisplaybreaks{\begin{align*}
J^{\mathcal F}_{\eta}(\chi_{A_q\cap B_q}) &= J^{\mathcal
F}_{\eta}(\chi_{A_q}\chi_{B_q}) \\ &= (1_{\mathcal E_s}\otimes
P_{\mathcal F^\bot_{q-s}}\otimes 1_{\mathcal E_{1-q}})(1_{\mathcal
E_q}\otimes P_{\mathcal F^\bot_{t-q}}\otimes 1_{\mathcal E_{1-t}}) \\
&= 1_{\mathcal E_s}\otimes P_{\mathcal F^\bot_{q-s}\otimes \mathcal
F^\bot_{t-q}}\otimes 1_{\mathcal E_{1-t}}. \end{align*}} %
For measurable sets $D_1,D_2, \cdots ,$ we claim that $$J^{\mathcal
F}_{\eta}(\chi_{\cup_i D_i})=\vee_i J^{\mathcal
F}_{\eta}(\chi_{D_i}).$$ Indeed $J^{\mathcal F}_{\eta}(\chi_{\cup_i
D_i})\geq J^{\mathcal F}_{\eta}(\chi_{D_i})$ for all $i,$ implying
$J^{\mathcal F}_{\eta}(\chi_{\cup_i D_i})\geq\vee_i J^{\mathcal
F}_{\eta}(\chi_{D_i}).$ For the reverse inequality, let
$D^\prime_1,D^\prime_2,\cdots $ be the sets constructed from
$D_1,D_2,\cdots $ by $D^\prime_1=D_1,$ $D^\prime_i=D_i\setminus
\cup_{k\leq {i-1}}D_k,$ $i\geq 2.$ Then $\cup_i D^\prime_i=\cup_i
D_i,$ $D^\prime_i \subset D_i$ and $D^\prime_1,D^\prime_2,\cdots $
are mutually disjoint. Then \allowdisplaybreaks{\begin{align*} J^{\mathcal F}_{\eta}(\chi_{\cup_i
D_i})&= J^{\mathcal F}_{\eta}(\chi_{\cup_i D^\prime_i}) \\ &=
J^{\mathcal F}_{\eta}(\sum_i \chi_{D^\prime_i}) \\ &= \vee_i
J^{\mathcal F}_{\eta}(\chi_{D^\prime_i}) \\ & \leq   \vee_i
J^{\mathcal F}_{\eta}(\chi_{D_i}).\end{align*}} %
So we have, \allowdisplaybreaks{\begin{align*} J^{\mathcal
F}_{\eta}(\chi_{ \{Z:\sharp\{Z\cap [s,t]\}\leq 1\}}) &= J^{\mathcal
F}_{\eta}(\chi_{(\cup_{q\in \mathbb Q\cap (s,t)}(A_q\cap B_q))^c}) \\
&= I_{\mathcal E_1}- \bigvee
_ {q\in \mathbb Q\cap (s,t)}J^{\mathcal F}_{\eta}(\chi_{A_q\cap B_q}) \\
&= I_{\mathcal E_1}- \bigvee_{q\in \mathbb Q\cap (s,t)} 1_{\mathcal
E_s}\otimes P_{\mathcal F^\bot_{q-s}\otimes \mathcal
F^\bot_{t-q}}\otimes 1_{\mathcal E_{1-t}} \\ &= I_{\mathcal E_1}-
1_{\mathcal E_s}\otimes P_{\vee_ {q\in \mathbb Q\cap (s,t)}(\mathcal
F^\bot_{q-s}\otimes \mathcal F^\bot_{t-q})}\otimes 1_{\mathcal
E_{1-t}} \\ &= 1_{\mathcal E_s}\otimes P_{[\vee_ {q\in \mathbb
Q\cap (s,t)}(\mathcal F^\bot_{q-s}\otimes \mathcal
F^\bot_{t-q})]^\bot}\otimes 1_{\mathcal E_{1-t}}. \end{align*}} %
We claim that
$$[\bigvee_{q\in \mathbb Q\cap (s,t)}[\mathcal F^\bot_{q-s}\otimes
\mathcal F^\bot_{t-q}]^\bot=[\bigvee_{s<r<t}\mathcal
F^\bot_{r-s}\otimes \mathcal F^\bot_{t-r}]^\bot.$$ For the moment
let us assume the claim. Then using the definition of $\mathcal
F^\prime,$ we get that \bea J^{\mathcal F}_{\eta}(\chi_{\{Z: Z\cap
[s,t]\leq 1\}})= P^{\mathcal F^\prime}_{s,t}. \eea For any partition
$\mathcal P=\{s=r_1<r_2<\cdots <r_k=t\}$ of $[s,t],$ define the set
$\mathcal A_{\mathcal P}=\{Z:Z\cap [r_i,r_{i+1}]\leq 1, 1\leq i \leq
k \}.$ Then the following identity holds: $$ \{Z:\sharp\{Z\cap
[s,t]\}< \infty\}=\cup_{\mathcal P}\mathcal A_{\mathcal P}.
$$ As $\check{\mathcal{F}}$ is the product system generated by the inclusion system $\mathcal F^\prime$, we have \allowdisplaybreaks{\begin{align*} P^{\check{\mathcal F}}_{s,t} &= \bigvee_{\mathcal P} \Pi _{i=1}^{n-1}P^{\mathcal F^\prime}_{r_i,r_{i+1}}\\
&= \bigvee_{\mathcal P} J^{\mathcal F}_{\eta}(\chi_{\mathcal A_{\mathcal P}})\\
&= J^{\mathcal F}_{\eta}(\chi_{ \cup_{\mathcal P} \mathcal A_{\mathcal P}})
\\ &=  J^{\mathcal F}_{\eta}(\chi_{\{Z:\sharp\{Z\cap [s,t]\}< \infty\}}).
\end{align*}} %
Now from the identity \bea \label{infty} \{Z:\check{Z}\cap
(s,t)=\emptyset\}=\{Z:\sharp\{Z\cap [s,t]\}<\infty\},\eea we get $$ J^{\mathcal F}_\eta(\chi_{\{Z:\check{Z}\cap [s,t]=\emptyset\}})= J^{\mathcal F}_\eta(\chi_{\{Z:\check{Z}\cap (s,t)=\emptyset\}})=P^{\check{\mathcal F}}_{s,t},{~}((s,t)\in I).$$ Therefore $$\mu^{\mathcal F}_{\eta}\{Z:\check{Z}\cap [s,t]=\emptyset \}=\eta(
P^{\check{\mathcal F}}_{s,t}) {~}{~},~([s,t]\subset I).
$$ Now it only remains to prove the claim. Clearly $$[\bigvee_{q\in \mathbb Q\cap (s,t)}\mathcal
F^\bot_{q-s}\otimes \mathcal F^\bot_{t-q}]^\bot \supset
[\bigvee_{s<r<t}\mathcal F^\bot_{r-s}\otimes \mathcal
F^\bot_{t-r}]^\bot.$$ Without loss of generality, we may assume
$0<s<t<1.$ We will use the continuity properties as in Proposition 3.18, \cite{Lie-random}. Let $$ P^0= 1_{\mathcal E_s}\otimes
P_{[\bigvee_{s<r<t}\mathcal F^\bot_{r-s}\otimes \mathcal
F^\bot_{t-r}]^\bot}\otimes 1_{\mathcal E_{1-t}}.$$  $$ Q^0=
1_{\mathcal E_s}\otimes P_{[\bigvee_{q\in \mathbb Q\cap
(s,t)}\mathcal F^\bot_{r-s}\otimes \mathcal
F^\bot_{t-r}]^\bot}\otimes 1_{\mathcal E_{1-t}}.$$ Clearly
$$P^0 \leq Q^0.$$ To prove the claim, it is enough to show that $P^0\geq Q^0.$ We have, \allowdisplaybreaks{\begin{align*} P_0 &= \bigwedge_{s<r<t}
[P^{\mathcal F}_{s,r}(1-P^{\mathcal F}_{r,t})+ (1-P^{\mathcal
F}_{s,r})P^{\mathcal F}_{r,t}+ P^{\mathcal F}_{s,t}]\\ &=
\bigwedge_{s<r<t} [P^{\mathcal F}_{s,r}+P^{\mathcal F}_{r,t}-
P^{\mathcal F}_{s,t}]. \end{align*}} %
Fix $x\in \mbox{range}~Q^0.$ Fix $r\in (s,t).$ Given $\epsilon>0,$ there is a $q\in \mathbb Q$ such that
$$\|P^{\mathcal F}_{s,r}x-P^{\mathcal F}_{s,q}x\|<\epsilon/2$$ and $$\|P^{\mathcal F}_{r,t}x-P^{\mathcal F}_{q,t}x\|<\epsilon/2 .$$
Now \allowdisplaybreaks{\begin{align*} \|x-[P^{\mathcal F}_{s,r}+P^{\mathcal F}_{r,t}- P^{\mathcal
F}_{s,t}]x\| &= \|[P^{\mathcal F}_{s,q}+P^{\mathcal F}_{q,t}-
P^{\mathcal F}_{s,t}]x-[P^{\mathcal F}_{s,r}+P^{\mathcal F}_{r,t}-
P^{\mathcal F}_{s,t}]x\|\\ & <  \epsilon/2+\epsilon/2 \\ & <
\epsilon. \end{align*}} %
So $$x \in \mbox{range}~[P^{\mathcal
F}_{s,r}+P^{\mathcal F}_{r,t}- P^{\mathcal F}_{s,t}]~{~} \mbox{for}
~\mbox{all}~ s<r<t.$$ This implies $$x\in \mbox{range}~P^0.$$ This
shows $P^0\geq Q^0$ and completes the proof. \qed

Suppose $(\mathcal E,W)$ is a product system and $u$ is a unit of
$(\mathcal E,W).$ For the product subsystem $\mathcal F_t=\mathbb Cu_t,$ we
wish to show that $\check{\mathcal F}$ is the type I part of $\mathcal E.$ To prove the result, we need the following lemmas.
\begin{lem}\label{X_s}
Suppose $(\mathcal E,W)$ is a product system and $u$ is a normalized
unit of $(\mathcal E,W).$ Then $u_s\otimes \mathcal F^\prime_t\subset
\mathcal F^\prime_{s+t}$ and $\mathcal F^\prime_s\otimes u_t\subset \mathcal F^\prime_{s+t}.$
\end{lem}
Proof: Suppose $x\in \mathcal F^\prime_t.$ consider the set
$$ A:=\{(z_1\otimes z_2):\langle z_1,u_r\rangle=0=\langle
z_2,u_{s+t-r}\rangle~,~ \mbox{for}~\mbox{some}~ r,~ 0<r<s+t \}.$$
Then we claim that $\overline{\mbox{span}}~A=\overline{\mbox{span}}~(A_1\cup A_2\cup A_3),$
where $$A_1=\{(y_1\otimes y_2\otimes y_3):\langle y_1,u_r\rangle=0 ,
\langle y_2,u_{s-r}\rangle\langle y_3,u_t\rangle=0,~
\mbox{for}~\mbox{some} ~0<r<s\},$$ $$A_2=\{(y_1\otimes y_2\otimes
y_3):\langle y_1,u_s\rangle\langle y_2,u_{r-s}\rangle=0 , \langle
y_3,u_{s+t-r}\rangle=0, \mbox{for}~\mbox{some} ~s<r<s+t\}$$ and $$ A_3=\{z_1\otimes z_2: \langle z_1,u_s\rangle=0, \langle z_2,u_t\rangle=0\}.$$ Suppose $y_1\otimes y_2\otimes y_3 \in A_1.$ That means for some $0<r<s,$ $\langle y_1,u_r\rangle=0,$ $\langle y_2,u_{s-r}\rangle\langle y_3,u_t\rangle=0.$ This implies $y_1 \in \mathcal E_r \ominus \mathbb Cu_r$ and $y_2\otimes y_3 \in \mathcal E_{s+t-r}\ominus \mathbb Cu_{s+t-r}.$ This shows $y_1\otimes y_2\otimes y_3 \in A.$ We obtain $A_1\subset A.$ Similarly, $A_2, A_3 \subset A.$  We obtain, $\overline{\mbox{span}}~A\supset \overline{\mbox{span}}~(A_1\cup A_2 \cup A_3).$ For the converse, let $z_1\otimes z_2\in A,$ with $\langle z_1,u_r\rangle=0, \langle z_2,u_{s+t-r}\rangle=0,$ $0<r<s.$ This implies $z_2 \in \overline{\mbox{span}} \{x_1\otimes x_2: x_1\in \mathcal E_{s-r}, x_2\in \mathcal E_t, \langle x_1\otimes x_2,u_{s+t-r}\rangle =0\}.$ Clearly $z_1\otimes x_1\otimes x_2 \in A_1.$ We get $z_1\otimes z_2 \in \overline{\mbox{span}}~A_1.$  Similarly, for $z_1\otimes z_2 \in A$ with $\langle z_1,u_r\rangle=0,$ $z_2,u_{s+t-r}\rangle=0,$ $s<r<s+t,$ we have $z_1\otimes z_2\subset \overline{\mbox{span}}~A_2.$  Therefore $\overline{\mbox{span}}~A\subset \overline{\mbox{span}}~(A_1\cup A_2\cup A_3).$ This proves the claim. Now suppose
$ y_1\otimes y_2\otimes y_3\in A_1$ be an arbitrary vector. Then
there is some $r_0,$ $0<r_0<s,$ such that $\langle
y_1,u_{r_0}\rangle=0 , \langle y_2,u_{s-r_0}\rangle\langle
y_3,u_t\rangle=0,$ and \allowdisplaybreaks{\begin{align*} \langle u_s\otimes x,y_1\otimes
y_2\otimes y_3\rangle &= \langle u_{r_0}\otimes u_{s-r_0}\otimes
x,y_1\otimes y_2\otimes y_3\rangle
\\ &= \langle
u_{r_0},y_1\rangle\langle u_{s-r_0},y_2\rangle\langle x,y_3\rangle \\
&= 0.\end{align*}} %
This shows that $u_s\otimes x \in {A_1}^\bot.$ Now let $y_1\otimes y_2\otimes y_3\in A_2$ be arbitrary. Then there is some $r_1,$ $s<r_1<s+t,$ such that $\langle y_1,u_s\rangle\langle
y_2,u_{r_1-s}\rangle=0 , \langle y_3,u_{s+t-r_1}\rangle=0.$ Now if
$\langle u_s,y_1\rangle=0,$ then the inner product $\langle
u_s\otimes x,y_1\otimes y_2\otimes y_3\rangle=0$ and if $\langle
u_s,y_1\rangle\neq 0,$ then $\langle y_2,u_{r_1-s}\rangle=0$ and
$\langle y_3,u_{s+t-r_1}\rangle=0.$ This is equivalent to
$y_2\otimes y_3 \in \tilde{F_t}.$ As $x\in F^\prime_t,$ the inner
product $\langle u_s\otimes x,y_1\otimes y_2\otimes y_3\rangle=0.$ For $z_1\otimes z_2\in A_3,$ it is easily seen that $\langle u_s\otimes x, z_1\otimes z_2 \rangle=0.$ Thus for arbitrary vector $z \in \overline{\mbox{span}}A,$ we have
$\langle u_s\otimes x, z\rangle=0.$ Hence $u_s\otimes x \in
F^\prime_{s+t}.$ Similarly $F^\prime_s\otimes u_t\subset
F^\prime_{s+t}.$ \qed

Define $X_t=F^\prime_t\ominus \mathbb Cu_t.$ From the previous lemma, it follows easily that
$u_s\otimes X_t\subset X_{s+t}$ and $X_s\otimes u_t\subset X_{s+t}.$
We identify the space $X_s$ as a subspace of $X_{s+t}$ by $x\mapsto
x\otimes u_t.$ This is an isometric embedding. Set
$X=\mbox{ind~limit}_{s>0} X_s.$ Denote the image of $x\in X_s,$ in
$X$ via $\overline{x}.$ For $t>0,$ define $S_t:X\rightarrow X$ via
$S_t(\overline{x})=\overline{u_t\otimes x},$ and set $S_0=id.$
\begin{lem}\label{L2}
Suppose $S_t$ is defined as above. Then $(S_t)_{t\geq 0}$ forms a
strongly continuous semigroup of isometries. Also
$(S_t)_{t\geq 0}$ is a pure semigroup i.e. $(S_t)_{t\geq 0}$ does
not have any unitary part.
\end{lem}
Proof: Clearly $(S_t)_{t\geq 0}$ is a semigroup of isometries. Now
to prove strong continuity of $t\mapsto S_t,$ it is enough to show
that for $x \in X_p, y \in X_q ,$ $0\leq t\leq 1,$
 $t \mapsto \langle \overline{x},S_t\overline{y}\rangle $ is continuous. Set $T>0$ such that $p<T,$ $q+t<T.$ Now \allowdisplaybreaks{\begin{align*} \langle \overline{x},S_t\overline{y}
\rangle &= \langle \overline{x},\overline{u_t\otimes y} \rangle \\
&= \langle x\otimes u_{T-p} , u_t \otimes y \otimes u_{T-q-t}
\rangle. \end{align*}} %
Set $z=x\otimes u_{T-p}$ and $w=y \otimes u_{T-q}.$ Then $z,w \in \mathcal E_T.$ Let $U^T_t$ in $\mathcal B(\mathcal E_T)$ be the unitary group $(U_t)_{t\in (0,T)}\subset
\mathcal B(\mathcal E_T)$ acting with regard to the representations
$\mathcal E_{T-t}\otimes \mathcal E_t\simeq\mathcal
E_T\simeq\mathcal E_t\otimes\mathcal E_{T-t}$ as flip: \bea \label{flip}
U^T_t(x_{T-t}\otimes y_t)=y_t\otimes x_{T-t}~,~(x_{T-t}\in \mathcal
E_{T-t} , y_t \in \mathcal \mathcal E_t).\eea Then $U^T_t(y\otimes
u_{T-q})=u_t\otimes y \otimes u_{T-q-t}.$ So we have,
$$ \langle \overline{x},S_t\overline{y} \rangle = \langle z , U^T_t
w\rangle.$$ Now an identical argument to the Proposition 3.11, \cite{Lie-random}, shows the map $t
\mapsto U^T_t$ is strongly continuous. Hence our result follows.

Now for the last part, it is equivalent to show that $\cap_{t\geq 0}
S_t(X)=\{0\}.$ we claim that $S_t(X_s)$ is orthogonal to $X_t$ for
every $s,t>0.$ Indeed, the claim follows from the fact that
$u_t\otimes X_s$ and $X_t\otimes u_s$ are orthogonal for every
$s,t>0.$ This implies that $S_t(X)$ is orthogonal to $X_t.$ Now if
$z\in \cap_{t\geq 0}S_t(X),$ we have $z \in {X_t}^\bot$ for every
$t>0$ $\Rightarrow$ $z \in X^\bot$ $\Rightarrow$ $z=0.$ \qed

\begin{thm}
Suppose $(\mathcal E,B)$ is a product system and $u$ is a normalized
unit of $(\mathcal E,B).$ Let $\mathcal F=(\mathcal F_t)_{t>0}$ be the product subsystem given by $\mathcal F_t:=\mathbb Cu_t.$ Then $(\check{\mathcal F},B)$ is the type I part
of $(\mathcal E,B).$
\end{thm}
Proof: From the definition of $\mathcal F^\prime_t,$ it is clear
that if $a$ is a root of the unit $u,$ then $a_t\in\check{\mathcal
F}_t.$ Now it follows from corollary \ref{unit-roots}, that
$(\check{\mathcal F},B)$ contains the type I part of $(\mathcal
E,B).$ On the other hand, we claim that
$$X_{s+t}=u_s\otimes X_t \oplus X_s \otimes u_t.$$ Indeed, from Lemma \ref{X_s}, we get $X_{s+t}\supset u_s\otimes X_t \oplus X_s \otimes u_t.$ For the reverse containment, we observe that,  $F^\prime_{s+t}\subset F^\prime_{s}\otimes F^\prime_t$ implies $X_{s+t}\subset (X_s \otimes u_t) \oplus (u_s\otimes X_t)\oplus (X_s\otimes X_t).$ So it is enough to show that $X_{s+t}\subset \mathcal E_{s+t}\ominus (X_s\otimes X_t).$ But this follows from the fact that $X_s\otimes X_t \subset \tilde{F}_{s+t}.$
So under the identification of $X_s$ inside $X_{s+t},$ we have for
every $s,t>0,$
$$X_{s+t}=X_s \oplus S_s(X_t).
$$ Taking limit as $t\uparrow \infty ,$ we get for every $s>0,$ $$X=X_s\oplus S_s(X).$$ Theorem 9.3, Chapter III, \cite{NaF-harmonic} states that every pure strongly continuous semigroup of isometries is unitarily equivalent to the unilateral shift semigroup on $L^2([0,t],K),$ for some Hilbert space $K.$ From Lemma \ref{L2}, there is a unitary $U:X\rightarrow
L^2([0,\infty),K)$ for some separable Hilbert space $K,$ such that
$US_tU^*=S^\prime_t,$ where $S^\prime_t$ is the unilateral shift
semigroup on $L^2([0,\infty],K).$  $$(S^\prime_t f)(x)= \left\{
\begin{array}{cc}
f(x-t) & \mbox{if}~x\geq t \\
0 & \mbox{otherwise}. \\
\end{array}
\right.$$ Now $L^2(0,\infty),K)$ decomposes
for every $s>0,$ as \allowdisplaybreaks{\begin{align*} L^2([0,\infty],K) &= U(X) \\ &=
U(X_s)\oplus US_s(X)
\\ &= U(X_s)\oplus US_sU^*U(X) \\ &= U(X_s)\oplus
S^\prime_s L^2([0,\infty],K) \\ &= U(X_s)\oplus L^2([s,\infty],K). \end{align*}} %
Therefore $U(X_s)=L^2([0,s],K).$ For each $t>0,$ Consider
the set $E_t= \mathbb C\Omega_t \oplus L^2([0,t],K)\subset
\Gamma_{sym}(L^2[0,t],K),$ where $\Omega_t=\Omega$ for all $t>0,$ is the vacuum vector.
Suppose $f\in L^2([0,s+t,K]).$ Let $g \in L^2([0,s],K)$ be defined
as $g=f|_{[0,s]}$ and $h\in L^2([0,t],K)$ be defined as
$$h(r)=f(r+s)~,~0\leq r \leq t.$$ Note that $S^\prime_sh=f|_{[s,s+t].}$ From the equation $f=g+S^\prime_sh,$ we have $W^{\Gamma}_{s,t}:\Gamma_{sym}(L^2[0,s+t],K)\rightarrow \Gamma_{sym}(L^2[0,s],K)\otimes \Gamma_{sym}(L^2[0,t],K)$ is given by $$W^{\Gamma}_{s,t}f=g\otimes \Omega_t + \Omega_s\otimes h. $$
For $\alpha \in \mathbb C$ and $f\in L^2([0,s+t],K),$
$$W_{s,t}^{\Gamma}|_{E_{s+t}}(\alpha\Omega_{s+t} \oplus f)=\alpha
(\Omega_s\otimes \Omega_t) \oplus (g\otimes\Omega_t + \Omega_s\otimes h)
\in E_s\otimes E_t. $$ So $(E_t,W^{\Gamma}_{s,t}|_{E_{s+t}})$ is an
inclusion system. Define $\Phi_t:\mathcal F^\prime_t \rightarrow
E_t$ by $$\Phi_t(\lambda u_t+x_t)=\lambda\Omega_s\oplus U|_{X_t}(x)$$
for $x\in X_t.$ We claim that $(\Phi_t)_{t>0}$ is an isometric
morphism of inclusion system. For $x\in X_{s+t},$ there are $y \in
X_t$ and $z\in X_s$ such that $W^{\mathcal E}_{s,t}x=u_s\otimes y+
z\otimes u_t.$ Under the identification on $X,$ we have $x=S_sy+z.$
So $Ux=US_sU^*Uy+Uz=S^\prime_sUy+Uz.$ Under the map
$W^{\Gamma}_{s,t}$ on $\Gamma_{sym}(L^2[0,s+t],K),$ we get
$S^\prime_tUy+Uz=W^{\Gamma
*}_{s,t}(\Omega_s \otimes Uy+Uz\otimes \Omega_t).$ This implies \allowdisplaybreaks{\begin{align*} W^{\Gamma *}_{s,t}(\Phi_s\otimes \Phi_t)W^{\mathcal
E}_{s,t}(\lambda u_{s+t}+x) &=  W^{\Gamma *}_{s,t}(\Phi_s\otimes
\Phi_t)(\lambda u_s\otimes u_t + u_s\otimes y+ z\otimes u_t) \\
&= W^{\Gamma
*}_{s,t} (\lambda (\Omega_s \otimes \Omega_t) \oplus (\Omega_s \otimes U|_{X_t}y + U|_{X_s}z\otimes
\Omega_t ))\\ &= \lambda \Omega_{s+t} \oplus S^\prime_sU|_{X_t}y+ U|_{X_s}z \\
&= \lambda \Omega_{s+t} \oplus US_sU^*U|_{X_t}y+ U|_{X_s}z \\ &= \lambda
\Omega_{s+t} \oplus U|_{X_{s+t}}(S_sy+z) \\ &= \lambda \Omega_{s+t} \oplus U|_{X_{s+t}}x\\
&= \Phi_{s+t}(\lambda u_{s+t}+x).\end{align*}} %
So $\mathcal F^\prime$ and $E$ are isomorphic as inclusion systems. So their generated product systems are isomorphic. As $E$ generates a type I product system, $\Gamma_{sym}(L^2[0,t],K),$ we have $\check{\mathcal F}$ generated by $F^\prime$ is a
type I product system of index $\mbox{dim}(K).$  \qed

Here we do a similar construction which generalize the cluster
construction.
 Suppose $\mathcal E$ is a product
system and $ F^1$ and $F^2$ are two inclusion subsystems of the
product system $\mathcal E.$ Consider for each $t>0,$ the space
$$G_t=\overline{\mbox{span}}\{x\otimes y: x\in \mathcal E_r\ominus F^1_r, y \in \mathcal E_{t-r}\ominus F^2_{t-r} ,
~\mbox{for}~\mbox{some}~0<r<t\}.$$ Define $G^\prime_t=\mathcal
E_t\ominus G_t.$
\begin{prop}
Let $G^\prime_t$ be defined as above. Then $G^\prime$ is an
inclusion system containing $F^1$ and $F^2.$
\end{prop}
Proof: First we will show that $F^1$ and $F^2$ is contained in
$G^\prime.$ Now fix $x\otimes y\in G_t.$ So $\langle
x,F^1_r\rangle=0,\langle y,F^2_{t-r}\rangle=0 ,$ for some $0<r<t.$
This implies that $\langle F^1_r\otimes F^1_{t-r},x\otimes
y\rangle=0,$ which with the fact $F^1_t\subset F^1_r\otimes
F^1_{t-r},$ for every $0<r<t$ proves that $F^1$ is contained in
$G^\prime.$ Similarly $ F^2\subset G^\prime.$ Suppose $x\in G^\prime_{s+t}.$ Now observe that $$G^\prime_s\otimes G^\prime_t=(\mathcal E_s\otimes G^\prime_t)\cap (G^\prime_s\otimes \mathcal E_t).$$  Now from the containment $F^1_{s+r}\subset F^1_s\otimes F^1_r\subset \mathcal E_s\otimes F^1_r,$ we  get $$\mathcal E_{s+r}\ominus F^1_{s+r}\supset \mathcal E_s\otimes (\mathcal E_r\ominus F^1_r).$$ For $y\in \mathcal E_s,$ $z_1\in \mathcal E_r\ominus F^1_r,$ $z_2\in \mathcal E_{t-r}\ominus F^2_{t-r},$ we get $y\otimes z_1 \in \mathcal E_{s+r}\ominus F^1_{s+r}.$ Consequently $\langle x,y\otimes z_1\otimes z_2 \rangle=0.$ We get $x\in \mathcal E_s\otimes G^\prime_t.$ Similarly, from the containment $F^2_{s+t-r^\prime}\subset F^2_{s-r^\prime}\otimes F^2_t\subset F^2_{s-r^\prime}\otimes \mathcal E_t,$ we get $$\mathcal E_{s+t-r}\ominus F^2_{s+t-r}\supset (\mathcal E_{s-r}\ominus F^2_{s-r})\otimes \mathcal E_t.$$ For $w_1\in \mathcal E_{r^\prime}\ominus F^1_{r^\prime},$ $w_2\in \mathcal E_{s-r^\prime}\ominus F^2_{s-r^\prime},$ $y^\prime \in \mathcal E_t,$ we get $w_2\otimes y^\prime \in \mathcal E_{s+t-r^\prime}\ominus F^2_{s+t-r^\prime}.$ Consequently $\langle x, w_1\otimes w_2\otimes y^\prime \rangle=0.$ We get $x\in G^\prime_s\otimes \mathcal E_t.$ Hence $x\in G^\prime_s\otimes G^\prime_t.$ Associativity of the inclusion system follows from the associativity of the product system. \qed

\begin{rem}
Note that if we take $F^1=F^2,$ then the product system generated by
$G$ is $\check{F^1}.$ Therefore it need not be the product system
generated by $F^1$ and $F^2.$
\end{rem}
\begin{center}
\textbf{Appendix A: More facts about inclusion systems}
\end{center}

Suppose $(E,\beta)$ is an inclusion system and $(\mathcal E,B)$ is its generated product system. We recall four basic properties of the inductive limit
construction.  (i) There exist canonical injections(isometries)
$i_\textbf{s}:E_\textbf{s}\rightarrow \mathcal E_t$  such that given
$\textbf{r}$ , $\textbf{s} \in J_t$ with $\textbf{r}\leq
\textbf{s}$ ,
$i_\textbf{s}\beta_{\textbf{s},\textbf{r}}=i_\textbf{r}$. ~(ii)
$\overline {\mbox{span}} \{i_\textbf{s}(a):a \in E_\textbf{s},
\textbf{s}\in J_t\}=\mathcal E_t.$ \label{closure} (iii) The
following universal property holds : Given a Hilbert space $\mathcal
G$ and isometries $g_\textbf{s}: E_\textbf{s}\rightarrow \mathcal G$
satisfying consistency condition
$g_\textbf{s}\beta_{\textbf{s},\textbf{r}}=g_\textbf{r}$ for all
$\textbf{r}\leq \textbf{s} \in J_t,$ there exists a unique isometry
$g:\mathcal E_t \rightarrow \mathcal G$ such that
$g_\textbf{s}=gi_\textbf{s}$ $\forall \textbf{s} \in J_t.$ (iv)
Suppose $K \subseteq J_t$ has the following property: Given
$\textbf{s} \in J_t,$ there exists $\textbf{t} \in K$ such that
$\textbf{s} \leq \textbf{t},$ then
$\mathcal E_t=\mbox{indlim}_{\textbf{r}\in K}E_\textbf{r}.$ Observe that $K$ is indeed a directed set with
the order inherited from $J_t.$ In other words,
$(x_{\textbf{s}})_{\textbf{s}\in K}$ is a subnet of
$(x_{\textbf{t}})_{\textbf{t}\in J_t}.$  We observe that the family of maps $i=(i_s)_{s>0}$ where $i_s:E_s\rightarrow \mathcal E_s$ satisfy the following: for $\textbf{s}=(s_1,\cdots,s_m)\leq \textbf{t}=(t_1,\cdots,t_n)\in J_t,$ we have  \bea \label{coproduct} B_{s,t}i_{\textbf{s}\smile\textbf{t}}=i_\textbf{s}\otimes i_\textbf{t}.\eea
 \bea \label{link} B_{\textbf{s},\textbf{t}}(i_{s_1}\otimes\cdots\otimes i_{s_m})= (i_{t_1}\otimes\cdots\otimes i_{t_n})\beta_{\textbf{s},\textbf{t}}.\eea

Here we prove the following important fact which we use repeatedly without reference.

\begin{lem}\label{subinclusion}
 Suppose $(\mathcal E,V)$ is a product system and $(F,\beta)$ is an inclusion subsystem of $(\mathcal E,V).$ Suppose $(\mathcal F,B)$ is the algebraic product system generated by $(F,\beta).$ Then $(\mathcal F,B)$ can be identified as a product subsystem of $(\mathcal E,V).$
 \end{lem}

 Proof: For every $s>0,$ $F_s$ is a closed subspace of $\mathcal E_s$ and for $s,t>0,$  $\beta_{s,t}=V_{s,t}|_{F_{s+t}}.$  Consider the family of isometries $(V^*_{t,\textbf{s}}|_{F_\textbf{s}}:F_\textbf{s}\rightarrow \mathcal E_t)_{\textbf{s}\in J_t}.$  Then for $\textbf{s}\leq \textbf{t}\in J_t,$ we have \allowdisplaybreaks{\begin{align*} V^*_{t,\textbf{t}}|_{F_\textbf{t}}\beta_{\textbf{s},\textbf{t}} & = V^*_{t,\textbf{t}}|_{F_\textbf{t}}V_{\textbf{s},\textbf{t}}|_{F_\textbf{s}} \\ &= V^*_{t,\textbf{t}}V_{\textbf{s},\textbf{t}}|_{F_\textbf{s}} \\ &=V^*_{t,\textbf{s}}|_{F_\textbf{s}}.  \end{align*} } %
 By the property (iii) listed above, for every $t>0,$ there is a unique isometry $j_t:\mathcal F_t\rightarrow \mathcal E_t$ such that for every $\textbf{s}\in J_t,$ $j_ti_\textbf{s}=V^*_{t,\textbf{s}}|_{F_\textbf{s}}$ where $i_\textbf{s}:F_\textbf{s}\rightarrow \mathcal F_s$ is the canonical inclusion. We claim that $j=(j_t)_{t>0}$ is an isometric morphism of algebraic product system  from $(\mathcal F,B)$ to $(\mathcal E,V).$ Indeed, for $\textbf{s}\in J_s$ and $\textbf{t}\in J_t,$ \allowdisplaybreaks{\begin{align*} (j_s\otimes j_t)(i_\textbf{s}\otimes i_\textbf{t}) &= j_\textbf{s}i_\textbf{s} \otimes j_\textbf{t}i_\textbf{t} \\ &= (V^*_{s,\textbf{s}}|_{E_\textbf{s}} \otimes V^*_{t,\textbf{t}}|_{E_\textbf{t}}) \\ &= V^*_{(s,t),\textbf{s}\smile \textbf{t}}|_{E_{\textbf{s}\smile\textbf{t}}} \\ &= V_{s,t} V^*_{s+t,\textbf{s}\smile\textbf{t}}|_{E_{\textbf{s}\smile\textbf{t}}} \\ &= V_{s,t}j_{s+t}i_{\textbf{s}\smile\textbf{t}}. \end{align*} %

 From Equation (\ref{coproduct}), we get $ (j_s\otimes j_t)B_{s,t}i_{\textbf{s}\smile \textbf{t}} = V_{s,t}j_{s+t} i_{\textbf{s}\smile \textbf{t}}.$
 As by property (iv) above, $$\overline{\mbox{span}~\{i_{\textbf{s}\smile\textbf{t}}(a):\textbf{s}\smile\textbf{t}\in J_s\smile J_t,a\in F_{\textbf{s}\smile\textbf{t}} \}}=\mathcal F_{s+t},$$ we get $$ (j_s\otimes j_t)B_{s,t} = V_{s,t}j_{s+t}.$$ This proves the claim. We may thus identify $\mathcal F$ as an algebraic product subsystem of $\mathcal E.$ Now suppose $U^{T,\mathcal E}=(U^{T,\mathcal E}_t)_{t\in \mathbb R}$ is the unitary group on $\mathcal E_T$ as in (\ref{flip}). Note that $U^{T,\mathcal F}_t=U^{T,\mathcal E}_t|_{\mathcal F_T}.$ From Proposition 3.11, [Lie] we get that $t\rightarrow U^{T,\mathcal E}_t$ is strongly continuous for every $T>0.$ Therefore being a restriction map, $t\rightarrow U^{T,\mathcal F}_t$ is also strongly continuous for every $T>0.$ Now an application of Theorem 7.7, \cite{Lie-random} shows that $(\mathcal F,B)$ is a product subsystem of $(\mathcal E,V).$

 \begin{center}
\textbf{Appendix B: Guichardte's picture of symmetric Fock space}
 \end{center}

 Suppose $K$ is a separable Hilbert space. Set $H=L^2(\mathbb R_,K),$ $H_t=L^2([0,t],K).$ We denote by $\Gamma_{sym}(H),$ $\Gamma_{sym}(H_t)$ the symmetric Fock spaces over $H$ and $H_t$ respectively. For $f\in \Gamma_{sym}(H_s)$ and $g\in \Gamma_{sym}(H_t),$ we define $W^\Gamma_{s,t}(e(f)\otimes e(g))=e(S_tf+g),$ where $S_t$ is defined as in Equation \eqref{S_t}. Then $\Gamma_{sym}(K):=(\Gamma_{sym}(H_t),W^\Gamma_{s,t})$ is a product system. Set $\Delta_n(t)$ is the set of all subsets of the interval $[0,t]$ of cardinality $n$ and $\Delta(t)$ is the set of all finite subsets of the interval $[0,t].$

 From the Lebesgue measure on the real line, we induce the Poisson measure $P$ on $\Delta(t)$ by $$P(E)=\delta_{\{\emptyset\}}(E)+\sum\limits_{1}^{\infty}\frac{1}{n!}\int\limits_{0}^{t}\cdots\int\limits_{0}^{t}dt_1\cdots dt_n \delta_{\{t_1,t_2,\cdots,t_n\}}(E).$$  Set $\mathcal F_t=L^2(\Delta(t),K,P)=\{g:\Delta(t)\rightarrow \Gamma(K): g(\sigma)\in K^{\otimes \#\sigma} , \int \|g\|^2dP < \infty \},$ where $\Gamma(K)$ is the full Fock space over $K$ and $\#\sigma$ is the cardinality of $\sigma.$  For $f\in L^2([0,t],K),$ denote $\hat{f}\in \mathcal F_t$ by $$\hat{f}(\sigma)=\left\{ \begin{array}{cc}
 1 & \mbox{if}~\sigma=\emptyset \\
 f(t_1)\otimes f(t_2)\otimes \cdots \otimes f(t_n) & \mbox{if}~ \sigma=\{t_1\leq t_2 \leq \cdots \leq t_n\}. \\
 \end{array}
 \right.$$ We claim that $\{\hat{f}:f\in L^2([0,t],K)\}$ is dense in $\mathcal F_t.$ Indeed, suppose $g\in \mathcal F_t$ and $\int\overline{\hat{f}}gdP=0,$ for all $f\in L^2([0,t],K).$ Then $\int\overline{\hat{1_{[a,b]}}}gdP=0,$ for all intervals $[a,b]\subset [0,t].$ It implies $\int\limits_{E}gdP=0$ for $E=\{\sigma \in \Delta(t):\sigma \subset [a,b]\}.$ As these sets are the cylinder sets for the sigma filed, we get $g=0.$
 Now under the map $e(f)\rightarrow \hat{f}$ we have the Hilbert space isomorphism $\Gamma_{sym}(H_t)\simeq \mathcal F_t.$ For $f\in \mathcal F_s,$ and $g\in \mathcal F_t,$ define $W^\mathcal F(f\otimes g) \in \mathcal F_{s+t}$ by $$W^\Gamma(f\otimes g)(\sigma)=f(\sigma\cap[t,s+t]-t)\otimes g(\sigma\cap [0,t]),~\sigma\in \Delta(s+t).$$ Then it is easily verified that the product systems $\Gamma_{sym}(K)=(\Gamma_{sym}(H_t),W^\Gamma)$ and $\mathcal F=(\mathcal F_t,W^\mathcal F)$ are isomorphic. Under this isomorphism, the vacuum vector $\Omega_t=e(0)$ is identified as $$\hat{0}(\sigma)=\left\{ \begin{array}{cc}
 1 & \mbox{if} ~\sigma=\emptyset \\
 0 & \mbox{otherwise}. \\
 \end{array}\right.$$

{\bf Acknowledgement.}  The third named author would like to thank DST-Inspire [IFA-13 MA-20] for financial support.

 \bibliography{ref}
 \bibliographystyle{amsplain}

\end{document}